\documentclass{elsarticle} % For numbered citations
\usepackage{amsmath,amssymb,amsthm,enumitem,mathtools}
\usepackage[hidelinks]{hyperref}
\usepackage{parskip}
\newtheorem{thm}{Theorem}[section]
\newtheorem{prop}[thm]{Proposition}
\newtheorem*{prop*}{Proposition}
\newtheorem{cor}[thm]{Corollary}
\newtheorem{lem}[thm]{Lemma}

\theoremstyle{definition}
\newtheorem{defn}[thm]{Definition}
\newtheorem{example}[thm]{Example}
\newtheorem{remark}[thm]{Remark}

\newcommand{\N}{\mathbb{N}}
\newcommand{\Z}{\mathbb{Z}}
\newcommand{\Q}{\mathbb{Q}}
\DeclareMathOperator{\Ap}{Ap}
\DeclareMathOperator{\Hom}{Hom}
\DeclareMathOperator{\Ker}{Ker}
\newcommand{\ff}{\mathfrak{f}}
\newcommand{\gapset}{\mathrm{H}}
\DeclareMathOperator{\Geom}{Geom}

\makeatletter
\NewCommandCopy\@@pmod\pmod
\DeclareRobustCommand{\pmod}{\@ifstar\@pmods\@@pmod}
\def\@pmods#1{\mkern4mu\left({\operator@font mod}\mkern 6mu#1\right)}
\makeatother
\begingroup
    \makeatletter
    \@for\theoremstyle:=definition,remark,plain\do{
        \expandafter\g@addto@macro\csname th@\theoremstyle\endcsname{
            \addtolength\thm@preskip\parskip
            }
        }
\endgroup

\begin{document}
\begin{frontmatter}
\title{Equidistribution Conditions for Gaps of Geometric Numerical Semigroups}
\author[1]{Caleb M. Shor\corref{cor1}}
\ead{caleb.shor@wne.edu}
\affiliation[1]{organization={Department of Mathematics, Western New England University},
    city={Springfield},
    state={MA},
    country={USA}}
\cortext[cor1]{Corresponding author}

\author[2]{Jae Hyung Sim}
\ead{simjhsim@bu.edu}

\affiliation[2]{organization={Department of Mathematics, Boston University},
    city={Boston},
    state={MA},
    country={USA}}

\begin{abstract}
    In 2008, Wang \& Wang showed that the set of gaps of a numerical semigroup generated by two coprime positive integers $a$ and $b$ is equidistributed modulo 2 precisely when $a$ and $b$ are both odd. Shor generalized this in 2022, showing that the set of gaps of such a numerical semigroup is equidistributed modulo $m$ when $a$ and $b$ are coprime to $m$ and at least one of them is 1 modulo $m$. In this paper, we further generalize these results by considering numerical semigroups generalized by geometric sequences of the form $a^k, a^{k-1}b, \dots, b^k$, aiming to determine when the corresponding set of gaps is equidistributed modulo $m$. With elementary methods, we are able to obtain a result for $k=2$ and all $m$. We then work with cyclotomic rings, using results about multiplicative independence of cyclotomic units to obtain results for all $k$ and infinitely many $m$. Finally, we take an approach with cyclotomic units and Dirichlet L-functions to obtain results for all $k$ and all $m$.
\end{abstract}
\end{frontmatter}
\section{Introduction}
Let $\N_0$ denote the set of non-negative integers, an additive monoid, and let $\N$ denote the set of positive integers. A \emph{numerical semigroup} $S$ is a submonoid of $\N_0$ with finite complement. Elements of the complement, which we denote $\gapset(S)$, are called the \emph{gaps} of $S$. The \emph{genus} of $S$, denoted $g(S)$, is the number of gaps of $S$. 

For any sequence of positive integers $g_1,\dots,g_\ell$, let $\langle g_1,\dots,g_\ell\rangle$ denote the additive monoid that they generate. It happens that $\langle g_1,\dots,g_\ell\rangle$ is a numerical semigroup if and only if $\gcd(g_1,\dots,g_\ell)=1$. Furthermore, every numerical semigroup $S$ arises this way. Given a numerical semigroup $S$, the minimal number of generating elements needed to generate $S$ is known as the embedding dimension of $S$, denoted $e(S)$.

In this paper, we investigate the problem of determining when the set of gaps of a numerical semigroup $S$ is \emph{equidistributed} modulo $m$, as was introduced in \cite{Shor22}. By this, we mean that the gap set $\gapset(S)$ contains the same number of elements in each congruence class modulo $m$.

As a starting point, consider any numerical semigroup $S$ with embedding dimension 2, so a numerical semigroup of the form $S=\langle a,b\rangle$ with $\gcd(a,b)=1$. A result of Wang \& Wang describes when $S$ has as many even gaps as odd gaps.
\begin{prop}[{\cite[Theorem 4.1]{WangWang2008}}]
    For $S=\langle a,b\rangle$ with $\gcd(a,b)=1$, the set $\gapset(S)$ is equidistributed modulo 2 if and only if $a$ and $b$ are both odd.
\end{prop}

We illustrate this proposition with two examples. The numerical semigroup $S_1=\langle 3,5\rangle$ has $\gapset(S_1)=\{1,2,4,7\}$, a set that is equidistributed modulo $2$. The numerical semigroup $S_2=\langle 4,5\rangle$ has $\gapset(S_2)=\{1,2,3,6,7,11\}$, a set that is not equidistributed modulo 2. However, this set has two elements in each congruence class modulo 3 and hence is equidistributed modulo 3. In \cite{Shor22}, this article's first author generalized the result of Wang \& Wang for embedding dimension 2 and any modulus $m$.
\begin{prop}[{\cite[Corollary 5.8]{Shor22}}]\label{prop:shor-equidistribution}
    Let $m\in\N$. For $S=\langle a,b\rangle$ with $\gcd(a,b)=1$, the set $\gapset(S)$ is equidistributed modulo $m$ if and only if $\gcd(ab,m)=1$ and at least one of $a$ and $b$ is 1 modulo $m$.
\end{prop}

Using the same examples $S_1$ and $S_2$, this result tells us that $\gapset(S_1)$ is equidistributed modulo 2 and modulo 4, and that $\gapset(S_2)$ is equidistributed modulo 3. The reader can verify these statements.

Thus, the question of determining when $\gapset(S)$ is equidistributed modulo $m$ is totally answered for numerical semigroups $S$ with embedding dimension 2. Quite often, going from embedding dimension 2 to embedding dimension greater than 2 can be quite a leap! For instance, for embedding dimension 2, the genus of $S=\langle a,b\rangle$ is $(a-1)(b-1)/2$, and the largest gap (known as the Frobenius number of $S$) is $ab-a-b$. If $S=\langle g_1,\dots,g_\ell\rangle$ for some $\ell>2$, then, as was shown by Curtis in \cite{Curtis1990}, there cannot exist a formula for the genus of $S$ as a polynomial in $g_1,\dots,g_\ell$. Furthermore, Ram\'irez-Alfons\'in showed, in \cite{RamirezAlfonsin96}, that computing the largest gap of $S$ (known as the Frobenius number of $S$) is NP-hard, under Turing reductions.

As a result, given any sequence of generating elements $g_1,\dots,g_\ell$ and a modulus $m$, we suspect it will be difficult to determine generally whether the set of gaps of $S=\langle g_1,\dots,g_\ell\rangle$ is equidistributed modulo $m$. In this paper, we generalize $S=\langle a,b\rangle$ by considering numerical semigroups generated by geometric sequences, which we refer to as \emph{geometric numerical semigroups}. 

\textbf{Question}. Consider the geometric numerical semigroup $S=\langle a^k, a^{k-1}b, \dots, b^k\rangle$, where $a,b,k\in\N$ with $\gcd(a,b)=1$. For which $m\in\N$ does it happen that $\gapset(S)$ is equidistributed modulo $m$? 

In Proposition~\ref{prop:shor-equidistribution}, a necessary condition for the set of gaps of $\langle a,b\rangle$ to be equidistributed modulo $m$ is that $\gcd(ab,m)=1$. We have the same necessary condition in the geometric case (Proposition~\ref{prop:gcd(ab,m)=1}). In order to answer our question, we begin by reframing the equidistribution condition for a geometric numerical semigroup as a polynomial congruence condition, working modulo $(x^m-1)$. In Proposition~\ref{prop:goal-reframed}, we find that, for $\gcd(ab,m)=1$, the set of gaps of $S=\langle a^k, a^{k-1}b, \dots, b^k\rangle$ is equidistributed modulo $m$ if and only if 
\[
    (x-1)\prod\limits_{i=1}^k\left(x^{ag_i}-1\right)\equiv\prod\limits_{i=0}^k\left(x^{g_i}-1\right)\pmod*{(x^m-1)},
\]
where $g_i=a^{k-i}b^i$ for $i=0,1,\dots,k$. With this, we can fully solve the $k=2$ case, concluding in Proposition~\ref{prop:result-k=2} that the set of gaps of $S=\langle a^2,ab,b^2\rangle$ is equidistributed modulo $m$ if and only if $\gcd(ab,m)=1$ and at least one of the following occurs: $a\equiv1\pmod*{m}$; $b\equiv1\pmod*{m}$; or $a^3\equiv b^3\equiv ab\equiv1\pmod*{m}$.

Unfortunately, the strategy that works for $k=2$ doesn't work for larger $k$. Instead, since we are working with polynomials modulo $(x^m-1)$, we move to work with cyclotomic fields, in particular $\Q(\zeta_d)$, where $\zeta_d$ is a primitive $d$th root of unity and $d$ is a divisor of $m$. The ring of integers of $\Q(\zeta_d)$ is $\Z[\zeta_d]$, and our primary objects of study in these rings are \emph{cyclotomic units} and dependence relations among them. It turns out that, for certain integers $m$, the set 
\[
    \left\{
        \frac{1-\zeta_m^h}{1-\zeta_m} : 1<h<m/2,\, \gcd(h,m)=1
    \right\}
\]
is (multiplicatively) independent. When this occurs (for instance, when $m$ is a power of a prime), we can fully answer the question about whether or not the set of gaps of a geometric numerical semigroup is equidistributed modulo $m$. (See Corollary~\ref{cor:prime-power-result} and what follows.) 

It happens that there are infinitely many $m$ for which this set is dependent. As a result, we need a more specialized approach to answer our question. In Section~\ref{sec:cyclotomic-units-and-L-values}, we make use of deep arithmetic properties of cyclotomic fields and the deep relationship between cyclotomic units and special values of Dirichlet L-functions (i.e., the base case of Stark's Conjectures) to resolve the equidistribution question for all $m$.

Our main result is Theorem~\ref{thm:main-result}, which says that the set of gaps of $S=\langle a^k, a^{k-1}b, \dots, b^k\rangle$ is equidistributed modulo $m$ if and only if any of seven conditions, which have to do with $m$, $k$, and the congruence classes of $a$ and $b$ modulo $m$, hold.

\subsection{Organization}
This paper is organized as follows. In Section~\ref{sec:prelim}, we state some preliminaries, showing how to use an \'Apery set to deduce an equivalent statement to equidistribution modulo $m$ for geometric semigroups in the form of a congruence of polynomials (with rational coefficients) modulo $(x^m-1)$ (Proposition~\ref{prop:goal-reframed}). 

In Section~\ref{sec:polynomial-congruence}, we demonstrate our first attempt at characterizing the congruence of polynomials as a list of congruences among powers of $a$ and $b$ modulo $m$. We fully characterize the case when $k=2$ (Proposition~\ref{prop:result-k=2}) in Section~\ref{subsec:equidistribution-k=2}. In Section~\ref{subsec:cyclo-dependence}, we reframe Proposition~\ref{prop:goal-reframed} to the study of dependence relations between cyclotomic units. This will provide enough leverage for a large array of moduli $m$ (Corollary~\ref{cor:prime-power-result} and {\cite[Theorem 3(2)]{Feng1982}}). However, as shown in Example~\ref{ex:cyclotomic-dependence}, we will require stronger tools to address the case of general $m$. 

In Section~\ref{sec:cyclotomic-units-and-L-values}, we will introduce these tools in the form of special values of Dirichlet L-functions and the representation theory of finite abelian groups. These tools will be used to prove Corollary~\ref{cor:more-general-result} which we will use in Section~\ref{sec:equidist-revisited} to fully characterize the equidistribution problem for geometric semigroups. Our main result is Theorem~\ref{thm:main-result}, which generalizes all previous results of this paper. 

Lastly, we provide a proof of Proposition~\ref{prop:gcd(ab,m)=1}  in Section~\ref{sec:gcd(ab,m)=1}. This result states that if the set of gaps of a geometric numerical semigroup $S$ is equidistributed modulo $m$, then $m$ is coprime to the generating elements of $S$.

\subsection{Questions for further study}
One generalization of a geometric sequence is that of a compound sequence, which is defined as follows: Let $a_1,\dots,a_k$ and $b_1,\dots,b_k$ be sequences of positive integers. Let $g_0=a_1a_2\cdots a_k$ and, for $1\le i\le k$, let $g_i = g_{i-1}\cdot b_i/a_i$. The sequence $g_0,\dots,g_k$ is called a compound sequence. If $\gcd(a_i,b_j)=1$ for all $i\ge j$, then $\gcd(g_0,\dots,g_k)=1$ and hence $S=\langle g_0,\dots,g_k\rangle$ is a numerical semigroup. (Compound sequences and the numerical semigroups that they generate first appeared in \cite{KiersONeillPonomarenko16}. See also \cite{GassertShor16}.) The question is then to determine all integers $m$ for the set of gaps of such a numerical semigroup $S$ is equidistributed modulo $m$.

Generalizing further, a compound sequence is a special case of what is known as a telescopic (or smooth) sequence. Numerical semigroups generated by telescopic sequences are called free numerical semigroups and they have appeared widely in the literature. (See \cite{Shor19a} and the references therein.) One can determine the moduli $m$ for which the sets of gaps of these numerical semigroups are equidistributed.

\subsection{Notations used in this paper}
For the convenience of the reader, we list here the notations for objects that we will use frequently throughout this paper.
\begin{itemize}
    \item $\N_0$: non-negative integers.
    \item $\N$: positive integers.
    \item $\Q$: rational numbers.
    \item $\overline{\Q}$: algebraic closure of $\Q$.
    \item $\gapset(S)$: gaps of $S$, i.e., $\N_0\setminus S$.
    \item $\langle g_1,\ldots,g_\ell\rangle$: $\N_0$-linear combinations of $g_1,\ldots,g_\ell$.
    \item $\Ap(S;t)$: (Definition~\ref{subsec:numerical-semigroups}) Ap\'ery set of $S$ relative to $t$.
    \item $\Geom(a,b;k)$: geometric sequence $a^k,a^{k-1}b,\ldots,ab^{k-1},b^k$.
    \item $C_m(x)$: the polynomial $(x^m-1)/(x-1)=1+x+\cdots+x^{m-1}=\sum_{i=0}^{m-1} x^i$.
    \item $\zeta_n$: primitive $n$-th root of unity.
    \item $\chi$: Dirichlet character.
    \item $\ff_\chi$: conductor of $\chi$.
    \item $\hat{G}$: character group of a group $G$.
    \item $\phi_{G/H}$: (Proposition~\ref{charproj}) averaging operator for functions on $G$.
    \item $U_m$: units modulo $m$, i.e., $(\Z/m\Z)^\times$.
    \item $\hat{U}_m^+$: $\chi\in \hat{U}_m$ such that $\chi(-1)=1$.
    \item $L(s,\chi)$: Dirichlet L-function associated to character $\chi$.
    \item $c_\chi$: (Proposition~\ref{Lvalue}) algebraic integer associated to $L(1,\chi)$.
    \item $V_m$: subspace of $\overline{\Q}$-valued functions on $U_m$ generated by nontrivial characters such that there exists a prime $p\mid m$ such that $p\nmid \ff_\chi$.
\end{itemize}

\section{Preliminaries}\label{sec:prelim}
We begin with two preliminary subsections. In the first, we describe a few tools for use with numerical semigroups and for geometric numerical semigroups in particular. (For a comprehensive treatment of numerical semigroups, see \cite{RosalesGarciaSanchez09}.) In the second, we translate the problem of detecting equidistribution modulo $m$ of the set of gaps of a geometric numerical semigroup to a polynomial congruence modulo $(x^m-1)$.

\subsection{Geometric numerical semigroups}\label{subsec:numerical-semigroups}
An important tool for use with any numerical semigroup is an Ap\'ery set of $S$.
\begin{defn}
    Let $S$ be a numerical semigroup and let $t$ be some nonzero element of $S$. The \emph{Ap\'ery set of $S$ relative to $t$}, denoted $\Ap(S;t)$, is defined as follows:
    \[\Ap(S;t) = \{s\in S : s-t\not\in S\}.\]
\end{defn}

In other words, $\Ap(S;t)$ consists of exactly $t$ elements; each is the smallest element of $S$ in its particular congruence class modulo $t$.

We have the following result which allows one to compute various quantities associated with a numerical semigroup $S$ given knowledge of an Ap\'ery set of $S$.
\begin{prop}\label{prop:tuenter-generalization}\cite[Theorem 2.3]{GassertShor17}
    For $S$ a numerical semigroup and $0\ne t\in S$, let $f$ be any arithmetic function. Then
    \[\sum\limits_{n\in \gapset(S)} \left(f(n+t)-f(n)\right) = \sum\limits_{n\in\Ap(S;t)}f(n) - \sum\limits_{i=0}^{t-1}f(i).\]
\end{prop}
As an example, the reader can confirm that using the function $f(n)=n$ with the above identity results in a formula for the genus of $S$.

We now focus on numerical semigroups generated by geometric sequences. Let $a,b,k\in\N$ with $\gcd(a,b)=1$. For $0\le i\le k$, let $g_i=a^{k-i}b^i$. Let $\Geom(a,b;k)$ denote the geometric sequence $g_0,g_1,\dots,g_k$. Then $S=\langle \Geom(a,b;k) \rangle$ is a \emph{geometric numerical semigroup}. If $a=1$ or $b=1$, then there are 0 gaps, so we conclude that the set of gaps is equidistributed modulo $m$ for all $m$. In what follows, we will assume that $a$ and $b$ are both greater than 1.

We begin with a few well-known results, starting with the genus of $S$. This was calculated in \cite[Proposition 7]{Shor05} and \cite[Theorem 6]{Shor11} in the context of the genus of an algebraic curve, as well as in \cite[Theorem 1(b)]{Tripathi08}.
\begin{prop}\label{prop:genus}
    For $\gcd(a,b)=1$, the genus of $S=\langle \Geom(a,b;k)\rangle$ is
    \begin{align*}
        g(S)
        &=\frac{(a-1)b^{k+1}-(b-1)a^{k+1}+(b-a)}{2(b-a)}\\
        &=\frac{(a-1)(a^{k-1}b+a^{k-2}b^2+\cdots+b^{k})-(a^k-1)}{2}.
    \end{align*}
\end{prop}
Next, we have an explicit description of the Ap\'ery set of $S$ relative to $g_0$, which we take from \cite[Proposition 3.6]{GassertShor17}. (See also \cite{Leher2007}.)
\begin{prop}\label{prop:apery-set}
    For $\gcd(a,b)=1$, let $S=\langle \Geom(a,b;k)\rangle$. The Ap\'ery set of $S$ relative to $g_0=a^k$ is 
    \[
        \Ap(S;g_0)=\left\{ \sum\limits_{i=1}^k g_in_i : 0\le n_1,\dots,n_k<a\right\}.
    \]
\end{prop}

\subsection{Polynomial congruence}
We translate the problem of determining equidistribution modulo $m$ to a congruence of two polynomials modulo $(x^m-1)$. The first step in doing so is the following observation.
\begin{prop}[{\cite[Proposition~2.15]{Shor22}}] \label{prop:set-equidist-iff}
    A finite set $T\subset\Z$ is equidistributed modulo $m$ if and only if 
    \[
        (x-1)\sum\limits_{t\in T}x^t\equiv0\pmod*{(x^m-1)}.
    \]
\end{prop}

\begin{remark}\label{rem:C_m(x)reduction}
    We can (and will) equivalently consider Proposition~\ref{prop:set-equidist-iff} as a congruence modulo $C_m(x):=\sum_{i=0}^{m-1} x^i$. This is justified because $(x-1)\sum_{t\in T}x^t$ is necessarily congruent to $0$ modulo $(x-1)$. Since $x^m-1=C_m(x)\cdot (x-1)$, by the Chinese Remainder Theorem, the congruence of Proposition~\ref{prop:set-equidist-iff} holds if and only if it holds modulo $C_m(x)$. In fact, since $(x-1)$ is invertible modulo $C_m(x)$, we can further claim that $T\subset \Z$ is equidistributed modulo $m$ if and only if $\sum_{t\in T}x^t\equiv 0\pmod*{C_m(x)}$.
\end{remark}

Via Propositions~\ref{prop:tuenter-generalization}, \ref{prop:apery-set}, and \ref{prop:set-equidist-iff}, we can reduce the problem of determining when the set of gaps of $S$ is equidistributed modulo $m$ to determining when a polynomial congruence holds modulo $(x^m-1)$. The final ingredient for the desired form of our congruence is the following proposition whose proof appears in the content of Section~\ref{sec:gcd(ab,m)=1}.
\begin{prop}\label{prop:gcd(ab,m)=1}
    Let $S=\langle \Geom(a,b;k)\rangle$. If $\gapset(S)$ is equidistributed modulo $m$, then $\gcd(ab,m)=1$.
\end{prop}

We now arrive at our main criterion for equidistribution modulo $m$.

\begin{prop}\label{prop:goal-reframed}
    For the geometric numerical semigroup $S=\langle \Geom(a,b;k)\rangle$ with $\gcd(a,b)=1$, the set $\gapset(S)$ is equidistributed modulo $m$ precisely when $\gcd(ab,m)=1$ and 
\begin{equation}\label{eqn:4}
    (x-1)\prod\limits_{i=1}^k\left(x^{ag_i}-1\right) \equiv \prod\limits_{i=0}^k \left(x^{g_i}-1\right)\pmod*{(x^m-1)}.
\end{equation}
\end{prop}
\begin{proof}
    By Proposition~\ref{prop:gcd(ab,m)=1}, if $\gcd(ab,m)>1$, then $\gapset(S)$ cannot be equidistributed modulo $m$. We proceed with the assumption that $\gcd(ab,m)=1$. We also note that there is nothing to check if $m=1$, so we assume $m>1$. Additionally, the two sides of Eq.~\eqref{eqn:4} are both 0 modulo $(x-1)$ and hence congruent to each other modulo $(x-1)$, so it will suffice to show that they are congruent modulo $C_m(x)$.
    
    Let $f:\N_0\to \Q[x]$ be the arithmetic function $f(n)=x^n$. Applying Proposition~\ref{prop:tuenter-generalization} with $f$ and $\Ap(S;g_0)$ as described in Proposition~\ref{prop:apery-set}, we get 
    \begin{equation}\label{eqn:1}
        \sum\limits_{n\in \gapset(S)}\left(x^{n+g_0}-x^n\right) = \sum\limits_{n_1=0}^{a-1}\sum\limits_{n_2=0}^{a-1}\cdots \sum\limits_{n_k=0}^{a-1} x^{g_1n_1+g_2n_2+\cdots+g_kn_k}-\sum\limits_{i=0}^{g_0-1}x^i.
    \end{equation}
    Cleaning this up, we equivalently find 
    \begin{equation}\label{eqn:2}
        (x^{g_0}-1)\sum\limits_{n\in \gapset(S)}x^n = \sum\limits_{n_1=0}^{a-1}x^{g_1n_1} \sum\limits_{n_2=0}^{a-1}x^{g_2n_2}\cdots \sum\limits_{n_k=0}^{a-1}x^{g_kn_k} - \sum\limits_{i=0}^{g_0-1} x^i.
    \end{equation}
    
    As in Remark~\ref{rem:C_m(x)reduction}, let $C_m(x)=\sum_{i=0}^{m-1}x^i\in\Q[x]$. Since $(x-1)$ has a multiplicative inverse modulo $C_m(x)$ and for each $i=0,1,\ldots,k$, the integers $g_i=a^{k-i}b^i$ and $ag_i$ are units modulo $m$, we can write 
    \[
        \sum_{n_i=0}^{a-1} x^{g_in_i}\equiv\frac{x^{ag_i}-1}{x^{g_i}-1}\pmod*{C_m(x)}.
    \]

    Applying this to Eq.~\eqref{eqn:2} modulo $C_m(x)$, we have
    \begin{equation}\label{eqn:3}
        (x^{g_0}-1)\sum\limits_{n\in \gapset(S)}x^n 
        \equiv \dfrac{x^{ag_1}-1}{x^{g_1}-1} \cdot \dfrac{x^{ag_2}-1}{x^{g_2}-1}\cdot \cdots\cdot \dfrac{x^{ag_k}-1}{x^{g_k}-1} - \dfrac{x^{g_0}-1}{x-1}\pmod*{C_m(x)}.
    \end{equation}

    We now consider the hypotheses. By Proposition~\ref{prop:set-equidist-iff} and Remark~\ref{rem:C_m(x)reduction}, $\gapset(S)$ being equidistributed modulo $m$ is equivalent to the congruence 
    \[
        \sum_{n\in\gapset(S)}x^n\equiv0\pmod*{C_m(x)}.
    \]
    Since $\gcd(g_0,m)=1$, this occurs if and only if the left side of Eq.~\eqref{eqn:3} is congruent to 0, which occurs if and only if the right side is congruent to 0. Finally, we conclude that the right side is congruent to 0 if and only if 
    \[
        (x-1)\prod\limits_{i=1}^k (x^{ag_i}-1) \equiv \prod\limits_{i=0}^k (x^{g_i}-1)\pmod*{C_m(x)},
    \]
    as desired.
\end{proof}
Proposition~\ref{prop:goal-reframed} provides a computational method for detecting equidistribution of $\gapset(S)$. We illustrate this with an example. 

\begin{example}\label{ex:checking-equidist-k=2}
    Suppose $a=9$, $b=11$, and $k=2$. Then the numerical semigroup $S=\langle 81,99,121\rangle$ has genus, by Proposition~\ref{prop:genus}, equal to 
    \[
        g(S)=\frac{(9-1)11^3-(11-1)9^3+(11-9)}{2(11-9)} = 840.
    \]
    If $m$ is an integer for which $\gapset(S)$ is equidistributed modulo $m$, then $m$ necessarily divides 840. To use Proposition~\ref{prop:goal-reframed}, we consider $m$ for which $\gcd(ab,m)=\gcd(99,m)=1$. In other words, we consider divisors of $840/{\gcd(840,99)}=280$.

    Let $L(x)$ and $R(x)$ denote (respectively) the left and right sides of Eq.~\eqref{eqn:4} from Proposition~\ref{prop:goal-reframed}. Then 
    \[
        L(x)=(x-1)(x^{891}-1)(x^{1089}-1)
        \quad\text{ and } \quad
        R(x) = (x^{81}-1)(x^{99}-1)(x^{121}-1).
    \] 
    For each candidate $m$, we can check to see if $L(x)\equiv R(x)\pmod*{(x^m-1)}$ by reducing exponents modulo $m$. We'll test $m=7$, $m=10$, and $m=70$. For each $m$, the reductions of $L(x)$ and $R(x)$ modulo $(x^m-1)$ appear in Figure~\ref{fig:reductions}.

    \begin{figure}
    \begin{center}
        \begin{tabular}{|c|c|c|}\hline
            $m$ & $L(x)$ modulo $(x^m-1)$ & $R(x)$ modulo $(x^m-1)$ \\ \hline 
            7 & $(x-1)(x^2-1)(x^4-1)$ & $(x^4-1)(x^1-1)(x^2-1)$ \\
            10 & $(x-1)(x^1-1)(x^9-1)$ & $(x^1-1)(x^9-1)(x^1-1) $\\
            70 & $(x-1)(x^{51}-1)(x^{39}-1)$ & $(x^{11}-1)(x^{29}-1)(x^{51}-1)$
            \\ \hline
        \end{tabular}
        \caption{Reductions of $L(x)$ and $R(x)$ modulo $(x^m-1)$, from Example~\ref{ex:checking-equidist-k=2}}
        \label{fig:reductions}
    \end{center}
    \end{figure}
    
    For $m=7$ and $m=10$, we have $L(x)\equiv R(x)\pmod*{(x^m-1)}$. Hence the set of gaps of $S$ is equidistributed modulo 7 and modulo 10.
    
    For $m=70$, it looks like the products are different modulo $(x^{70}-1)$. We expand out and compute 
    \[
        L(x)-R(x)
        \equiv x^1+x^{10}-x^{11}-x^{20}-x^{29}+x^{39}-x^{52}+x^{62}\pmod*{(x^{70}-1)},
    \] 
    which is nonzero modulo $(x^{70}-1)$. Thus, $L(x)\not\equiv R(x)\pmod*{(x^{70}-1)}$. We conclude that the set of gaps of $S$ is not equidistributed modulo 70. (This allows us to further conclude that the set of gaps of $S$ is not equidistributed modulo any multiple of 70.)
\end{example}

\section{Conditions for a particular polynomial congruence}\label{sec:polynomial-congruence}
In this section, we address the following question which arises naturally in light of Example~\ref{ex:checking-equidist-k=2}.

Given two integer sequences $c_1,\dots,c_n$ and $d_1,\dots,d_n$ and a positive integer $m$ for which 
\begin{equation}\label{eqn:prod-to-exponents-congruence}
    \prod\limits_{i=1}^n \left(x^{c_i}-1\right)\equiv \prod\limits_{i=1}^n\left(x^{d_i}-1\right)\pmod*{(x^m-1)},
\end{equation}
what can we say about the two integer sequences modulo $m$? 

Just as in Remark~\ref{rem:C_m(x)reduction}, we note here that the congruence in Eq.~\eqref{eqn:prod-to-exponents-congruence} holds modulo $(x^m-1)$ if and only if it holds modulo $C_m(x)$. In our work below, it will sometimes be easier to work modulo $C_m(x)$ because certain polynomials (like $(x-1)$) are units modulo $C_m(x)$ and hence we can use their multiplicative inverses to simplify computations.

\subsection{Results when the exponents are units}\label{subsec:exp-are-units}

To start, observe that any permutation of the sequence $c_1,c_2,\dots,c_n$ will leave the left side of Eq.~\eqref{eqn:prod-to-exponents-congruence} unchanged. That is, for any permutation $\sigma$ in the symmetric group $S_n$, 
\[
    \prod\limits_{i=1}^n\left(x^{c_i}-1\right) = \prod\limits_{i=1}^n\left(x^{c_{\sigma(i)}}-1\right).
\] 
We can think of the sequences $c_1,\dots,c_n$ and $c_{\sigma(1)},\dots,c_{\sigma(n)}$ as effectively the same for our purposes. This leads us to the define two sequences to be congruent modulo $m$.

\begin{defn}[Congruent sequences modulo $m$]
    We say \emph{two integer sequences $c_1,\dots,c_n$ and $d_1,\dots,d_n$ are congruent modulo $m$} if there is some permutation $\sigma$ in the symmetric group $S_n$ such that $c_i\equiv d_{\sigma(i)}\pmod*{m}$ for all $i=1,2,\dots,n$.
\end{defn}

We will now consider Eq.~\eqref{eqn:prod-to-exponents-congruence} and see what we can conclude. To start, suppose one of the $c_i$ terms is 0 modulo $m$. Then the left side of Eq.~\eqref{eqn:prod-to-exponents-congruence} is 0 modulo $(x^m-1)$, which implies the right side is as well, and hence some $d_j$ is 0 modulo $m$. As both sides of the polynomial congruence are 0 modulo $(x^m-1)$, any relations between the remaining terms in the two sequences are completed washed away. To avoid this sort of situation, we will assume that the terms in our sequences are nonzero modulo $m$.

We make one more assumption. In our main criterion for equidistribution (Proposition~\ref{prop:goal-reframed}), all of the exponents in Eq.~\eqref{eqn:4} are coprime to $m$. We will make that simplifying assumption here, assuming that $\gcd(c_i,m)=\gcd(d_i,m)=1$ for all $i=1,2,\dots,n$. In particular, this means that $(x^{c_i}-1)$ and $(x^{d_i}-1)$ are units modulo $(x^m-1)$.

In the following proposition, we show that if Eq.~\eqref{eqn:prod-to-exponents-congruence} holds and $n\le 3$, then the two sequences must be congruent modulo $m$. It turns out that for $n\ge4$, this is not the case. (See Example~\ref{ex:n=4}.) Looking ahead, we will show in Section~\ref{sec:cyclotomic-units-and-L-values} that, for any $n$, if the polynomial congruence holds then the two sequences must be PM-congruent modulo $m$ (which we define in the next subsection).

\begin{prop}\label{prop:prods-match-exponents-match}
    Let $m\in\N$ and $n\in \{1,2,3\}$. For any two sequences of integers $c_1,\dots,c_n$ and $d_1,\dots,d_n$ relatively prime to $m$,
    Eq.~\eqref{eqn:prod-to-exponents-congruence} holds if and only if the sequences $c_1,\dots,c_n$ and $d_1,\dots,d_n$ are congruent modulo $m$.
\end{prop}
\begin{proof}
    The reverse direction is clearly true for all $n$. For the forward direction, any two units are congruent to each other if $m=1$ or $m=2$, so we will assume $m\ge 3$. This also means all $c_i$ and $d_j$ are nonzero modulo $m$.

    When $n=1$, Eq.~\eqref{eqn:prod-to-exponents-congruence} implies $x^{c_1}\equiv x^{d_1}\pmod*{(x^m-1)}$, so $c_1\equiv d_1\pmod*{m}$. This concludes the $n=1$ case.
    
    When $n=2$, Eq.~\eqref{eqn:prod-to-exponents-congruence} implies
    \begin{equation}\label{eqn:proof-n=2}
        (x^{c_1}-1)(x^{c_2}-1)
        \equiv 
        (x^{d_1}-1)(x^{d_2}-1)
        \pmod*{(x^m-1)}.
    \end{equation}
    Expanding and shuffling some terms yields
    \[
        x^{c_1+c_2}+x^{d_1}+x^{d_2}
        \equiv 
        x^{d_1+d_2}+x^{c_1}+x^{c_2}
        \pmod*{(x^m-1)},
    \]
    and hence the sequences $d_1+d_2,c_1,c_2$ and $c_1+c_2,d_1,d_2$ are congruent modulo $m$. Since $c_1$ and $c_2$ are nonzero modulo $m$, $c_1,c_2\not\equiv c_1+c_2\pmod*{m}$. Thus, $c_1+c_2\equiv d_1+d_2\pmod*{m}$ which implies that the sequences $c_1,c_2$ and $d_1,d_2$ are congruent modulo $m$. This concludes the $n=2$ case.

    When $n=3$. Eq.~\eqref{eqn:prod-to-exponents-congruence} is of the form
    \begin{equation}\label{eqn:proof-n=3}
        \prod\limits_{i=1}^3\left(x^{c_i}-1\right)
        \equiv
        \prod\limits_{i=1}^3\left(x^{d_i}-1\right)
        \pmod*{C_m(x)}.
    \end{equation}
    If it happens that $c_i\equiv d_j\pmod*{m}$ for some $i,j\in\{1,2,3\}$, then we can cancel $(x^{c_i}-1)$ and $(x^{d_j}-1)$ from each side modulo $C_m(x)$ and apply the $n=2$ result to the resulting two-term products to conclude the sequences $c_1,c_2,c_3$ and $d_1,d_2,d_3$ are congruent modulo $m$. It will therefore suffice to show that $c_i\equiv d_j\pmod*{m}$ for some $i,j$. We proceed by assuming for contradiction that $c_i\not\equiv d_j\pmod*{m}$ for all choices of $i,j\in \{1,2,3\}$.

    Now, consider Eq.~\eqref{eqn:proof-n=3} modulo $(x^m-1)$. Expanding and moving terms around, the following two polynomials are congruent modulo $(x^m-1)$:
    \[x^{c_1+c_2+c_3}+x^{d_1+d_2}+x^{d_1+d_3}+x^{d_2+d_3}+x^{c_1}+x^{c_2}+x^{c_3},\]
    \[x^{d_1+d_2+d_3}+x^{c_1+c_2}+x^{c_1+c_3}+x^{c_2+c_3}+x^{d_1}+x^{d_2}+x^{d_3}.\]
    Hence, the 7-term sequences $d_1+d_2+d_3,c_1+c_2,c_1+c_3,c_2+c_3,d_1,d_2,d_3$ and $c_1+c_2+c_3,d_1+d_2,d_1+d_3,d_2+d_3,c_1,c_2,c_3$ are congruent modulo $m$. 
    
    By our assumption for contradiction, we observe that $c_1,c_2,c_3$ is congruent modulo $m$ to a length-3 subsequence of $d_1+d_2+d_3,c_1+c_2,c_1+c_3,c_2+c_3$. This implies we have $c_i\equiv c_j+c_k\pmod*{m}$ and $c_j\equiv c_i+c_k\pmod*{m}$ for some distinct choices of $i,j,k\in\{1,2,3\}$. Summing these congruences together yields $2c_k\equiv 0\pmod*{m}$, which contradicts the earlier assumption that $c_k$ is a unit modulo $m\ge 3$. Thus, some $c_i$ is congruent to some $d_j$ modulo $m$, as desired. This concludes the $n=3$ case.
\end{proof}

\subsection{Equidistribution of the gaps of a geometric numerical semigroup for \texorpdfstring{$k=2$}{k=2}}\label{subsec:equidistribution-k=2}
We will now apply Proposition~\ref{prop:prods-match-exponents-match} to determine when $\gapset(S)$ is equidistributed modulo $m$ for the case where $k=2$; i.e., for $S=\langle a^2,ab,b^2\rangle$.

\begin{prop}\label{prop:result-k=2}
    For $a,b\in\N$ with $a,b>1$ and $\gcd(a,b)=1$, let $S=\langle a^2,ab,b^2\rangle$. Let $m\in\N$. Then $\gapset(S)$ is equidistributed modulo $m$ if and only if $\gcd(ab,m)=1$ and either $a\equiv1\pmod*{m}$ or $b\equiv1\pmod*{m}$ or $a^3\equiv b^3\equiv ab\equiv1\pmod*{m}$.
\end{prop}
\begin{proof}
    ($\implies$) If $\gcd(ab,m)>1$, then by Proposition~\ref{prop:gcd(ab,m)=1}, $\gapset(S)$ is not equidistributed modulo $m$. Thus, $\gcd(ab,m)=1$.
    
    By Proposition~\ref{prop:goal-reframed}, $\gapset(S)$ is equidistributed modulo $m$ if and only if 
    \begin{equation}\label{eqn:k=2-case}
        (x-1)(x^{a^2b}-1)(x^{ab^2}-1)\equiv (x^{a^2}-1)(x^{ab}-1)(x^{b^2}-1)\pmod*{(x^m-1)}.
    \end{equation}
    By Proposition~\ref{prop:prods-match-exponents-match}, since these exponents are all units modulo $m$, this congruence holds if and only if the three-term sequences $1,a^2b,ab^2$ and $a^2,ab,b^2$ are congruent modulo $m$. From here, we have three possibilities for the term that $ab$ is congruent to.

    First, assume $ab\equiv1\pmod*{m}$. If $a^2\equiv a^2b\pmod*{m}$, then $b\equiv1\pmod*{m}$. If instead $a^2\equiv ab^2\pmod*{m}$, then $a^3\equiv a^2b^2\equiv (1)^2\equiv 1\pmod*{m}$ and $b^3\equiv a^3b^3\equiv (ab)^3\equiv 1^3\equiv1\pmod*{m}$.
    
    Alternatively, if $ab\equiv a^2b\pmod*{m}$, then $a\equiv1\pmod*{m}$. By symmetry, we also have that if $ab\equiv ab^2\pmod*{m}$, then $b\equiv1\pmod*{m}$.

    Collecting the results, for all three possible congruences for $ab$, we have $a\equiv1\pmod*{m}$; or $b\equiv1\pmod*{m}$; or $a^3\equiv b^3\equiv ab\equiv1\pmod*{m}$.

    ($\impliedby$) Suppose $\gcd(ab,m)=1$. Then in order to show that $\gapset(S)$ is equidistributed modulo $m$, it is enough to show that Eq.~\eqref{eqn:k=2-case} holds.
    
    If $a\equiv1\pmod*{m}$ or $b\equiv1\pmod*{m}$, then Eq.~\eqref{eqn:k=2-case} holds and we are done. If $a^3\equiv b^3\equiv ab\equiv 1\pmod*{m}$, then the left side becomes
    \[
        (x-1)(x^{a^2b}-1)(x^{ab^2}-1)\equiv (x-1)(x^a-1)(x^b-1)\pmod*{(x^m-1)}.
    \] 
    Since $b\equiv a^{-1}\equiv a^2\pmod*{m}$ and $a\equiv b^{-1}\equiv b^2\pmod*{m}$, the right side becomes 
    \[
        (x^{a^2}-1)(x^{ab}-1)(x^{b^2}-1)\equiv (x^b-1)(x^1-1)(x^a-1)\pmod*{(x^m-1)}.
    \]
    so Eq.~\eqref{eqn:k=2-case} holds, as desired.
\end{proof}

Note that the third condition in Proposition~\ref{prop:result-k=2} may be simplified slightly. If $ab\equiv1\pmod*{m}$, then $a^3\equiv1\pmod*{m}$ if and only if $b^3\equiv1\pmod*{m}$. Thus, it is enough to check that $a^3\equiv ab\equiv1\pmod*{m}$.

\begin{example}\label{ex:checking-equidist-k=2-again}
    In Example~\ref{ex:checking-equidist-k=2}, with $a=9$, $b=11$, and $k=2$, we checked to see if the set of gaps of $S=\langle 9^2, 9\cdot11, 11^2\rangle$ is equidistributed modulo $m$ for $m=7$, $m=10$, and $m=70$. (Since we need $\gcd(ab,m)=1$, we only considered divisors of $g(S)/{\gcd(ab,g(S))}=280$.)
    
    We can now use Proposition~\ref{prop:result-k=2} to streamline these calculations. We want divisors $m$ of 280 such that one of the following occurs: $a\equiv 9\equiv1\pmod*{m}$; $b\equiv 11\equiv1\pmod*{m}$; or $a^3\equiv 9^3\equiv1\pmod*{m}$ and $ab\equiv 99\equiv1\pmod*{m}$. The first two cases imply $m$ divides 8 or 10. The third case says $m$ divides $\gcd(9^3-1,99-1)=14$. Thus, $\gapset(S)$ is equidistributed modulo $m$ if and only if $m$ divides 8, 10, or 14 which is consistent with the conclusion of Example~\ref{ex:checking-equidist-k=2}.
\end{example}

Unfortunately, the approach for Proposition~\ref{prop:result-k=2} does not allow us to say anything for $k\ge3$. The issue comes from the fact that the analog of Proposition~\ref{prop:prods-match-exponents-match} for $n\ge4$ is not true, as we illustrate with the following example.

\begin{example}\label{ex:n=4}
    For $m=4$, writing $x^{-1}\equiv x^3\pmod*{(x^4-1)}$, we have 
    \[(x^1-1)^4\equiv (-x^1)^4 (-1+x^{-1})^4\equiv (-1)^4 x^{4}(x^3-1)^4\equiv 1\cdot1\cdot(x^3-1)^4\pmod*{(x^4-1)}\]
    but $1,1,1,1$ and $3,3,3,3$ are not congruent sequences modulo 4.
\end{example}

While the two sequences in Example~\ref{ex:n=4} are not congruent modulo 4, they are still closely related, as the terms in the two sequences are additive inverses of each other. This motivates a generalization of our notion of congruent sequences.

\begin{defn}[PM-congruent sequences modulo $m$]
    We say \emph{two integer sequences $c_1,\dots,c_n$ and $d_1,\dots,d_n$ are plus-minus-congruent, or PM-congruent, modulo $m$} if there is some permutation $\sigma$ in the symmetric group $S_n$ such that, for each $i=1,2,\dots,n$, either $c_i\equiv d_{\sigma(i)}\pmod*{m}$ or $c_i\equiv -d_{\sigma(i)}\pmod*{m}$.
\end{defn}

Back to our question at the start of this section, certainly if two sequences are congruent modulo $m$ then the congruence in Eq.~\eqref{eqn:prod-to-exponents-congruence} holds. And there are pairs of PM-congruent sequences, such as in Example~\ref{ex:n=4}, for which the congruence holds. Our goal moving forward is to show that if Eq.~\eqref{eqn:prod-to-exponents-congruence} holds, then the two sequences are PM-congruent modulo $m$. We will require a new set of tools to tackle the general case, making some progress in the next subsection and ultimately finishing it off in Section~\ref{sec:cyclotomic-units-and-L-values}.

\subsection{Dependence relations in cyclotomic rings}\label{subsec:cyclo-dependence}

We now investigate conclusions drawn from Eq.~\eqref{eqn:prod-to-exponents-congruence} using cyclotomic units. Our notations follow \cite[Chapter 8]{washington}. Let $m$ be a positive integer such that $m\not\equiv2\pmod*{4}$, let $\zeta_m$ be a primitive $m$th root of unity with $K_m=\Q(\zeta_m)$, the $m$th cyclotomic field. The ring of integers of $K_m$ will be denoted $\mathcal{O}_m=\mathcal{O}_{K_m}=\Z[\zeta_m]$. 

To motivate our focus on cyclotomic units, we first notice that $x^m-1=\prod_{\ell\mid m}\Phi_\ell(x)$ where $\Phi_\ell$ is the minimal polynomial of $\zeta_\ell$ for each $\ell$. The Chinese Remainder Theorem then converts Eq.~\eqref{eqn:prod-to-exponents-congruence} to the following equality in $K_\ell$ for each $\ell\mid m$ where $\ell>1$:
\[
    \prod\limits_{i=1}^n \Big(\zeta_\ell^{c_i}-1\Big)= \prod\limits_{i=1}^n\Big(\zeta_\ell^{d_i}-1\Big).
\]
In what follows, we will see that if $m$ is a prime power, the properties of elements of the form $\zeta_\ell^i-1$ where $1\le i\le \ell-1$ will dictate when Eq.~\eqref{eqn:prod-to-exponents-congruence} is satisfied.

Letting $V_m$ be the multiplicative group generated by 
\[
    \{\pm\zeta_m,\, 1-\zeta_m^a : 1\le a\le n-1\}
\]
and $E_m$ be the group of units of $\Q(\zeta_m)$, the group of \emph{cyclotomic units} of $\Q(\zeta_m)$ is defined as $C_m=V_m\cap E_m$. As we will see below, there are infinitely many integers $m$ for which the cyclotomic units generate a subgroup of finite index of the full unit group of $\mathcal{O}_m$. When this occurs, we can fully determine the possible pairs of sequences $c_1,\dots,c_n$ and $d_1,\dots,d_n$, all units modulo $m$, for which Eq.~\eqref{eqn:prod-to-exponents-congruence} holds. 

By Dirichlet's Unit Theorem, the rank of $\mathcal{O}_m^\times$ (i.e., the unit group of $\mathcal{O}_m$) is $\varphi(m)/2-1$. The cases when $\mathcal{O}_m^\times/(\mathcal{O}_m\cap C_m)$ is finite occurs exactly when the rank of the set 
\[
    T_m=\left\{\frac{1-\zeta_m^h}{1-\zeta_m} : 1< h\le m-1\right\},
\]
is $\varphi(m)/2-1$. (Note that $T_m\subset\mathcal{O}_m$.) We can see two types of relations between elements of $T_m$. First, $1-\zeta_m^a = -\zeta_m^a(1-\zeta_{m}^{m-a})$ for all $a$. Second, suppose $\gcd(l,m)=d$. Then we can write $l=dl'$ and $m=dm'$ for coprime integers $l'$ and $m'$. Furthermore, there is some unit $u$ modulo $m$ such that $l\equiv du\pmod*{m}$. Following Bass' approach from \cite{Bass1966}, given that 
\[
    t^d-1
    =\prod\limits_{i=1}^d\Big(t-\zeta_d^{i}\Big)
    =\prod\limits_{i=1}^d\Big(t-\zeta_m^{m'i}\Big),
\]
if we set $t=\zeta_m^{-1}$ and then multiply both sides through by $\zeta_m^{d}$, we find 
\[1-\zeta_m^d=\prod\limits_{i=1}^d\left(1-\zeta_n^{m'i+1}\right).\]
Since $u$ is a unit modulo $m$, we replace $\zeta_m$ with $\zeta_m^u$ to conclude
\[
    1-\zeta_m^l = \prod\limits_{i=1}^d\left(1-\zeta_m^{(m'i+1)u}\right).
\] 
If $d>1$ then we have just obtained a nontrivial dependence relation.

As a result, instead of taking $1< h\le m-1$ in $T_m$, we may restrict values of $h$ to be units modulo $m$ from the first half of the interval from $1$ to $m-1$. Let 
\[
    T_m' = \left\{\frac{1-\zeta_m^h}{1-\zeta_m} : 1<h<m/2,\, \gcd(h,m)=1 \right\}.
\]

The following result follows from \cite[Theorem 1]{Ramachandra1966}. 
\begin{prop}\label{prop:prime-power-prod-to-exp}
    If $m$ is a power of a prime, then $T_m'$ is an independent set.
\end{prop}

We can now prove the prime power case.

\begin{cor}\label{cor:prime-power-result}
    Suppose $m$ is a power of a prime. If $c_1,\dots,c_n$ and $d_1,\dots,d_n$ are two sequences of integers coprime to $m$ such that 
    \[
        \prod\limits_{i=1}^n\left(x^{c_i}-1\right)
        \equiv
        \prod\limits_{i=1}^n\left(x^{d_i}-1\right)
        \pmod*{(x^m-1)},
    \]
    then the sequences $c_1,\dots,c_n$ and $d_1,\dots,d_n$ are PM-congruent modulo $m$.
\end{cor}
\begin{proof}
    For any integer $x$ coprime to $m$, let $\tilde{x}$ be the unique integer such that $1\le \tilde{x}<m/2$ such that either $x\equiv \pm\tilde{x}\pmod*{m}$. As in the beginning of this section, if 
    \[
        \prod\limits_{i=1}^n\left(x^{c_i}-1\right)
        \equiv
        \prod\limits_{i=1}^n\left(x^{d_i}-1\right)
        \pmod*{C_m(x)},
    \]
    then we have the following equality in $K_{m}$:
    \begin{equation}\label{eqn:proof-primepower}
        \prod\limits_{i=1}^n\left(\zeta_{m}^{c_i}-1\right)
        = \prod\limits_{i=1}^n\left(\zeta_{m}^{d_i}-1\right).
    \end{equation}
    Since $1-\zeta_m^{-h}=-\zeta_m^{-h}(1-\zeta_m^h)$, we have the following equality up to a factor of a root of unity.
    \[
        \prod\limits_{i=1}^n\frac{1-\zeta_{m}^{\tilde{c}_i}}{1-\zeta_m}
        = \prod\limits_{i=1}^n\frac{1-\zeta_{m}^{\tilde{d}_i}}{1-\zeta_m}.
    \]
    By Proposition~\ref{prop:prime-power-prod-to-exp}, this equality implies that for each integer $1<h<m/2$ coprime to $m$, 
    \[
        \#\{i : \tilde{c}_i\equiv h\pmod*{m}\}=\#\{i : \tilde{d}_i\equiv h\pmod*{m}\}.
    \]
    Lastly, we note that 
    \begin{align*}
        \#\{i : \tilde{c}_i\equiv 1\pmod*{m}\}&=n-\sum_{1<h<m/2}\#\{i : \tilde{c}_i\equiv h\pmod*{m}\}\\
        &=n-\sum_{1<h<m/2}\#\{i : \tilde{d}_i\equiv h\pmod*{m}\}\\
        &=\#\{i : \tilde{d}_i\equiv 1\pmod*{m}\}
    \end{align*}
    which implies that $c_1,\ldots,c_n$ and $d_1,\ldots,d_n$ are PM-congruent modulo $m$.
\end{proof}

At this point, our problem is solved for the case where $m$ is a prime power. We could now repeat a similar argument as in the proof of Proposition~\ref{prop:result-k=2} to deduce the characterizing properties for $\gapset(S)$ to be equidistributed modulo a prime power $m$. However, since the argument is not significantly different for the composite case given the correct analog of Corollary~\ref{cor:prime-power-result}, we defer until Section~\ref{sec:equidist-revisited} for this argument.

Instead we will now revisit the set $T_m'$ defined above in the case where $m$ is not a prime power. As we will see, there are infinitely many non-prime power integers $m$ for which $T_m'$ is an independent set (and hence we obtain a result similar to Corollary~\ref{cor:prime-power-result}). 

Relations between cyclotomic units have appeared in the literature quite often throughout the 20th century. (See \cite{Bass1966}, \cite{Ramachandra1966}, \cite{Ennola1972}, \cite{Feng1982}, and others.) We refer to the work of Feng \cite{Feng1982} in particular, where we find explicit criteria on $m$ for which the set $T_m'$ is independent. (Note that Feng states this theorem for non-prime power $m$. We note that it holds for prime power $m$ as well, so we have not made that distinction in our statement of the theorem.)

\begin{thm}[{\cite[Theorem 3(2)]{Feng1982}}]
    Suppose $m$ is a positive integer with $m\not\equiv2\pmod*{4}$, and suppose $m$ has prime factorization $m=p_1^{a_1}p_2^{a_2}\cdots p_l^{a_l}$ for distinct primes $p_1,p_2,\dots,p_l$, and positive exponents $a_1,a_2,\dots,a_l$. The set $T_m'$ is an independent set if and only if the unit group $U_{m/p_i^{a_i}}$ is generated multiplicatively by $-1$ and $p_i$ for all $i=1,2,\dots,l$.
\end{thm}

Immediately afterward (in \cite[Theorem 4]{Feng1982}), Feng shows that $T_m'$ is independent for non-prime power $m$ if and only if $m$ belongs to one of five families of integers. In this case, we obtain an analogue of Corollary~\ref{cor:prime-power-result}. Among positive integers $m$ up to 120 (still with $m\not\equiv2\pmod*{4}$), $T_m'$ is independent for all $m$ except: 39, 55, 56, 65, 68, 84, 91, 112, 117, 120. (These integers appear in sequence \href{https://oeis.org/A387182}{A387182} in the On-Line Encyclopedia of Integer Sequences \cite{oeis}.)

\begin{example}\label{ex:cyclotomic-dependence}
    We will show that $T_{39}'$ is a dependent set. For $\zeta_{39}$ a primitive 39th root of unity, let $G$ be the subgroup of $U_{39}$ generated by 35. Then $G=\{35, 16, 14, 22, 29, 1\}$. The reader can verify that 
    \begin{equation}\label{eqn:dep-reln-n=39}
        \prod\limits_{g\in G}\left(1-\zeta_{39}^g\right) 
        = 1.
    \end{equation}
    Since 35, 22, and 29 are greater than 39/2, we factor out $(-\zeta_{39}^{35})(-\zeta_{39}^{22})(-\zeta_{39}^{29})=-\zeta_{39}^{8}$ and then raise both sides to the 78th power to obtain
    \[
        \prod\limits_{g\in G'}\left(1-\zeta_{39}^g\right)^{78}
        =1
    \]
    for $G'=\{-35,16,14,-22,-29,1\}=\{4,16,14,17,10,1\}\subset U_{39}$.
    
    Since $\zeta_{39}^2$ is also a primitive 39th root of unity, we can replace $\zeta_{39}$ with $\zeta_{39}^2$ in Eq.~\eqref{eqn:dep-reln-n=39} to find
    \[
        \prod\limits_{h\in H}\left(1-\zeta_{39}^{h}\right)
        =1
    \]
    where $H=2G=\{31,32,28,5,19,2\}$. Following the same process as above, we find
    \[
        \prod\limits_{h\in H'}\left(1-\zeta_{39}^h\right)^{78}
        =1
    \]
    for $H'=\{-31,-32,-28,5,19,2\}=\{8,7,11,5,19,2\}\subset U_{39}$.

    We now have a dependence relation among elements of $T_{39}'$:
    \[
        \prod\limits_{g\in G', g\ne 1} \left(\frac{1-\zeta_{39}^g}{1-\zeta_{39}}\right)^{78}
        \prod\limits_{h\in H'} \left(\frac{1-\zeta_{39}^h}{1-\zeta_{39}}\right)^{-78}
        =1.
    \]
\end{example}

Unfortunately for us, there are infinitely many $m$ for which $T_m'$ is dependent. (In particular, $T_m'$ is dependent if $m$ has at least four prime factors.) At this point, we are unable to get an analogue of Corollary~\ref{cor:prime-power-result} for all $m$. We need a more robust set of numbers to replace the role of cyclotomic units. These numbers will be special values of Dirichlet L-functions which will be the main object of study in Section~\ref{sec:cyclotomic-units-and-L-values}.

\section{Cyclotomic Units and L-values}\label{sec:cyclotomic-units-and-L-values}

The goal of this section will be to prove Lemma~\ref{mainlemma} to 
uniformly address the conditions of Eq.~\eqref{eqn:prod-to-exponents-congruence}
for all $k$. 

To motivate Lemma~\ref{mainlemma} in the context of our problem, we first 
reinterpret the congruence of Eq.~\eqref{eqn:prod-to-exponents-congruence}
as the vanishing of a function. Taking the sequences $(c_i)$ and $(d_i)$ in
Eq.~\eqref{eqn:prod-to-exponents-congruence}, let $C,D:U_m\to\Z$ be the functions
\begin{align*}
    C(z)&=\#\{1\le i\le n : c_i\equiv z\pmod*{m}\}\\
    D(z)&=\#\{1\le i\le n : d_i\equiv z\pmod*{m}\}.
\end{align*}
Then, Eq.~\eqref{eqn:prod-to-exponents-congruence} is equivalent
to the congruence
\[
    \prod_{z\in U_m}(x^z-1)^{C(z)-D(z)}\equiv 1\pmod*{C_m(x)}.
\]
Lemma~\ref{mainlemma} will allow us to show that if $C(z)-D(z)$ 
satisfies the above congruence and some other mild conditions,
then $C(z)-D(z)=0$. This will allows us to conclude that $C(z)=D(z)$
which implies the congruence of the sequences $c_1,\ldots,c_n$ and $d_1,\ldots,d_n$ modulo $m$.

The proof of Lemma~\ref{mainlemma} relies on the general fact that 
for any finite group $G$, the set of characters of $G$ forms a 
$\overline{\Q}$-basis for the space of $\overline{\Q}$-valued class 
functions on $G$. Since all $\overline{\Q}$-valued functions are class 
functions if $G$ is abelian, we can decompose any arbitrary function 
$f:U_m\to \Z$ as a sum $f= \sum_{\chi}\alpha_{f,\chi}\cdot \chi$ where 
the sum is over all characters of $U_m$, i.e., homomorphisms $U_m\to 
\overline{\Q}^\times$. Showing that $f=0$ is then equivalent to showing
that $\alpha_{f,\chi}=0$ for each $\chi$.

We now outline the proof of Lemma~\ref{mainlemma}.
$\Q[x]$ being a PID, $\Q[x]/(C_m)$ is isomorphic to the product 
$\bigoplus_{d\mid m, d\neq 1}\Q(\zeta_d)$ where $\zeta_d$ is a primitive
$d$-th root of unity. By the Chinese Remainder Theorem, 
Eq.~\eqref{eqn:prod-to-exponents-congruence} is equivalent to a collection 
of equalities of the form $\prod_{a\in U_m}(\zeta_d^a-1)^{f(a)}=1$ for 
each $d\mid m$.

The elements $\zeta_d-1$ are called cyclotomic units and enjoy a rich history 
in the arithmetic of 
cyclotomic fields. In particular, they are the first examples of Stark 
units, which means they arise in the value of Dirichlet L-functions at $s=1$.
Explicitly, one can embed $\Q(\zeta_d)$ into $\mathbb{C}$ by mapping 
$\zeta_d$ to $e^{2\pi i/d}$ and take the complex logarithm of 
$\prod_{a\in U_m} (\zeta_d^a-1)^{f(a)}=1$ to yield the equality
\[
\sum_{\chi} \alpha_{f,\chi}\cdot c_\chi \cdot E_{\chi,d}\cdot 
L(1,\overline{\chi}) =0\] 
where $E_{\chi,d},c_\chi\in \overline{\Q}$ are related to the
value of Dirichlet characters and $c_\chi\neq 0$.

A deep theorem (Theorem~\ref{Lindependence}), whose proof involves relations 
between cyclotomic units and the transcendence of their logarithms, informs
us that for the subset of characters $\chi$ such that $E_{\chi,m}\neq 0$, the 
$\overline{\Q}$-independence of $\chi$ is preserved by the map 
$\chi\mapsto L(1,\overline{\chi})$. Thus, we can project $f$ to this
subspace of characters which allows us to deduce that $\alpha_{f,\chi}=0$ 
for such $\chi$. The proof is then an exercise in carefully considering 
the correct divisor of $m$ to focus on for each $\chi$. 

An important remark is that we require the full force of the congruence of 
polynomials in our assumption; the statement of Theorem~\ref{polyequality} is 
simply not true if we replace the congruence with just the equality 
$\prod_{a\in U_m}(\zeta_d^a-1)^{f(a)}=1$ for a particular $d\mid m$, even 
if $d=m$. For example, consider the subgroup $H\subset U_5\times U_{13}=U_{65}$
generated by all of $U_5$ and the subgroup of $U_{13}$ generated by the 
residue class of $5$. Then, $U_{65}/H$ is a group of size $3$. Furthermore, 
the product $\prod_{a\in H}(\zeta_{65}^a-1)$ is equal to $1$. This means
that if $f:U_{65}\to \Z$ is equal to $1$ on $H$, $-1$ on a nontrivial coset of 
$H$, and $0$ on the remaining coset, we have
\[
\prod_{a\in U_{65}} (\zeta_{65}^a-1)^{f(a)} = 1.
\]
However, $f$ is clearly not equal to $0$ by construction. In the scope of 
our proof, this is because there exists two nontrivial characters $\chi$
which factor through $U_{65}/H$ such that $\alpha_{f,\chi}\neq 0$.
We can see this fact by observing that
\[
\prod_{a\in U_{65}} (\zeta_{13}^a-1)^{f(a)}\neq 1.
\]

\subsection{Representation theory of finite abelian groups}
Let $G$ be a finite abelian group which decomposes as 
$\bigoplus_{i=1}^n \Z/p_i^{e_i}\Z$ for some collection of primes $p_i$ and 
$e_i\in \N$. Let $\hat{G}$ be the set of characters 
$\hat{G}=\Hom(G,\overline{\Q}^\times)$. 
\begin{prop}\label{charproj}
    Let $H\subset G$ be a normal subgroup and let $\phi_{G/H}$ be the operator on 
    $\overline{\Q}$-valued functions on $G$ which for any $f:G\to \overline{\Q}$ and 
    $x\in G$
    \[
    \phi_{G/H}(f)(x)=\frac{1}{\# H}\sum_{h\in H} f(hx).
    \]
    Then, for all $\chi\in \hat{G}$, 
    \[
    \phi_{G/H}(\chi)=\begin{cases}
        \chi & \text{if }\chi \text{ factors through the quotient }G/H\\
        0 & \text{otherwise.}
    \end{cases}
    \]
\end{prop}

\begin{proof}
    If $\chi\in \hat{G}$ factors through the quotient $G/H$, then for all $g\in H$
    and $x\in G$, we have $\chi(x)=\chi(gx)$. Thus, $\phi_{G/H}(\chi)=\chi$.

    If $\chi\in \hat{G}$ does not factor through $G/H$, then $\chi$ is nontrivial when
    restricted to $H$. Thus, the sum $\sum_{h\in H}\chi(h)=0$. It follows that
    for all $x\in G$, \[
    \phi_{G/H}(\chi)(x)=\sum_{h\in H}\chi(hx)=\sum_{h\in H}\chi(h)\chi(x)
    =\chi(x)\sum_{h\in H}\chi(h)=0.
    \]
\end{proof}

If we have a surjective homomorphism of groups $\psi:G_1\to G_2$, we will often write
$\phi_{G_2}$ for the operator $\phi_{G_1/\Ker(\psi)}$.

\begin{prop}\label{cycproj}
    Let $f:U_m\to \Z$ be a function and $d\mid m$ a positive divisor not equal to $1$.
    Then,
    \[
    \prod_{a\in U_m}(\zeta_d^a-1)^{f(a)} = \prod_{b\in U_d}(\zeta_d^a-1)^{\tfrac{\# U_m}{\# U_d}(\phi_{U_{m/d}}f)(b')},
    \]
    where $b'$ denotes any choice of element in $U_m$ congruent to $b$ in $U_d$.
\end{prop}

\begin{proof}
    Let $H\subset U_m$ be the kernel of the projection $U_m\to U_{m/d}$. Then, for all 
    $a\in U_m$, $\zeta_d^a-1=\zeta_d^{ah}-1$. This gives us for all $b\in U_d$,
    \begin{align*}
        \prod_{a\in \psi^{-1}(b)}(\zeta_d^a-1)^{f(a)}
        &= (\zeta_d^b-1)^{\sum_{a\in \psi^{-1}(b)}f(a)}\\
        &= (\zeta_d^b-1)^{\sum_{h\in H}f(b'\cdot h)}\\
        &= (\zeta_d^b-1)^{\#H\cdot (\phi_{U_{m/d}}f)(b')}.
    \end{align*}
    Taking the product over all $b\in U_d$ then gives us the desired equality.
\end{proof}

\subsection{Dirichlet L-functions}
For this section, let $m>1$ be a positive integer that is not congruent to $2$ modulo $m$ and let $\zeta_m$ be a primitive $m$th root of unity.

\begin{defn}
    A \emph{Dirichlet character} is a function $\chi:\Z\to \mathbb{C}$ such that there exists an
    integer $\ff_\chi$, called the \emph{conductor}, such that \begin{itemize}
        \item for all $n\in \Z$ such that $\gcd(n,\ff_\chi)\neq 1$, $\chi(n)=0$,
        \item for all $a,b\in \Z$ such that $\gcd(ab,\ff_\chi)=1$, $\chi(ab)=\chi(a)\cdot \chi(b)$,
        \item for all $a\in \Z$, $\chi(a)=\chi(a+\ff_\chi)$.
    \end{itemize} 

    We denote by $\hat{U}_m^+$ the set of nontrivial Dirichlet characters $\chi$ such that 
    $\ff_\chi\mid m$ and $\chi(-1)=1$.
\end{defn} 

We make the observation that characters of $U_m$ are canonically in one-to-one correspondence
with Dirichlet characters whose conductor divides $m$ by reducing mod $m$/extending by $0$ to $\Z$.
We liberally switch between these frameworks as needed.

\begin{thm}\label{Lindependence} (Corollary~25.6 of \cite{murtyrath})
    The set $\{L(1,\chi) : \chi\in \hat{U}_m^+\}$ is linearly independent 
    over $\overline{\Q}$.
\end{thm}

\begin{prop}(pg 149 of \cite{washington}) \label{Lvalue}
    Let $\chi$ be an even Dirichlet character such that $\mathfrak{f}_\chi\mid m$. Then, there exists
    a nonzero algebraic number $c_\chi\in\overline{\Q}$ such that
    \[
    \sum_{a\in U_m} \chi(a)\log |1-e^{2\pi i a/m}|
    =c_\chi \cdot L(1,\overline{\chi})\cdot \prod_{p\mid m,\text{ prime}}(1-\chi(p)).
    \]
\end{prop}

\begin{defn}
    For all $\chi\in \hat{U}_m$ and $d\mid m$, let $E_{\chi,d}$ be the algebraic number 
    $\prod_{p\mid d,\text{ prime}}(1-\chi(p))$. 
    
    Define $V_m$ as the subspace of $\overline{\Q}$-valued functions on $U_m$ generated by 
    \textit{nontrivial} characters $\chi \in \hat{U}_m$ such that there exists a prime 
    $p\mid m$ not dividing $\ff_\chi$. Define $W_m$ as the subspace of $\overline{\Q}$-valued
    functions on $U_m$ generated by \textit{nontrivial} characters $\chi\in \hat{U}_m^+$
    such that for all primes $p\mid m$, $p\mid \ff_\chi$.
\end{defn}

\begin{prop}\label{eulervanish}
    Let $\chi$ be a nontrivial Dirichlet character with $\ff_\chi\mid m$. If 
    $E_{\chi,m}=0$, then $\chi\in V_m$.
\end{prop}

\begin{proof}
    If $E_{\chi,m}=0$, there exists a prime $p\mid m$ such that $1-\chi(p)=0$.
    In particular, $\chi(p)$ is nonzero, so $p$ does not divide $\mathfrak{f}_\chi$.
\end{proof}

\begin{lem}\label{mainlemma}
    Let $f:U_m\to \Z$ be an even function satisfying $\sum_{a\in U_m}f(a)=0$ whose
    decomposition into characters is of the form $f=\sum_{\chi\in \hat{U}_m}\alpha_\chi \cdot \chi$.
    Then, the following two statements are true:
    \begin{enumerate}
        \item If $\chi$ is an odd character or $\chi$ is the trivial character,
        then $\alpha_\chi=0$.
        \item If $\sum_{a\in U_m} f(a)\cdot \log |1-e^{2\pi ia/m}| = 0$, 
        then $f\in V_m$.
    \end{enumerate}
\end{lem}

\begin{proof}
    For the first statement, we notice that $\phi_{\{1\}}f = \tfrac{1}{\# U_m}\sum_{a\in U_m}f(a)=0$, 
    so $\alpha_{\mathbf{1}}=0$ where $\mathbf{1}$ is the trivial character. Furthermore, 
    $\phi_{U_m/\{\pm 1\}}f = f$ since $f$ is even. Since $\phi_{U_m/\{\pm 1\}}$ is a linear operator
    that kills off all odd characters, $\alpha_{\chi}$ must have been $0$ for all odd $\chi$. 

    For the second statement, we use Proposition~\ref{Lvalue} and the first statement to notice that 
    \begin{align*}
        0&=\sum_{a\in U_m} f(a)\cdot \log |1-e^{2\pi ia/m}|\\
        &= \sum_{a\in U_m}\sum_{\chi\in \hat{U}_m^+}
        \alpha_\chi \cdot \chi(a) \cdot \log |1-e^{2\pi i a/m}|\\
        &= \sum_{\chi\in \hat{U}_m^+} \alpha_\chi\cdot c_\chi\cdot E_{\chi,m}\cdot L(1,\overline{\chi}).
    \end{align*}
    Theorem~\ref{Lindependence} implies that for all $\chi\in \hat{U}_m^+$, $\alpha_\chi\cdot c_\chi\cdot 
    E_{\chi,m}=0$. Since $c_\chi\neq 0$ for all $\chi\in \hat{U}_m^+$, $\alpha_\chi=0$ or $E_{\chi,m}=0$. 
    By Proposition~\ref{eulervanish}, if $\chi\in \hat{U}_m^+\cap W_m$, then $\alpha_\chi=0$. Thus, 
    $\alpha_\chi\neq 0$ only if $\chi\in \hat{U}_m^+\cap V_m$, so $f\in V_m$.
\end{proof}

\begin{thm}\label{polyequality}
    Suppose $f,g:U_m\to \Z$ are even functions such that $\sum_{h\in U_m}f(h)=
    \sum_{h\in U_m}g(h)$. If \[
    \prod_{h\in U_m}(x^h-1)^{f(h)}\equiv \prod_{h\in U_m}(x^h-1)^{g(h)}\pmod*{\tfrac{x^m-1}{x-1}},
    \]
    then $f=g$.
\end{thm}

\begin{proof}
    Since $\Q[x]$ is a principal ideal domain, the given congruence is equivalent
    to the statement that for all positive $d\mid m$ where $d\neq 1$, we have 
    \[\prod_{a\in U_m} (\zeta_d^a-1)^{f(a)}=\prod_{a\in U-m}(\zeta_d^a-1)^{g(a)}\]
    in $\Q(\zeta_d)$. 

    Consider the function $f-g:U_m\to \Z$ and its unique decomposition 
    $\sum_{\chi\in \hat{U}_m}\alpha_\chi\cdot \chi$. We now deduce that $\alpha_\chi=0$
    for all $\chi\in \hat{U}_m$.

    Let $\chi\in \hat{U}_m$. If $\chi$ is odd or trivial, Lemma~\ref{mainlemma}(1)
    tells us that $\alpha_\chi=0$ since by assumption, $f-g$ must be even and 
    $\sum_{a\in U_m}(f-g)(a)=\sum_{a\in U_m}f(a)-g(a)=0$.

    If $\chi$ is even, let $d\mid m$ such that $\gcd(d,m/d)=1$ and if $p\mid d$, then
    $p\mid \ff_\chi$. By assumption, we have the equality
    \[
    \prod_{a\in U_m}(\zeta_{d}^a-1)^{(f-g)(a)} = 1.
    \]
    By Proposition~\ref{cycproj}, this is equal to the product
    \[
    \prod_{b\in U_d}(\zeta_d^b-1)^{F(b)} = 1,
    \]
    where $F:U_d\to \Q$ is the function on $U_d$ that lifts to $\phi_{U_d}(f-g)$.
    Since $d$ was chosen so that $\gcd(d,m/d)=1$, $F$ is still an even function.
    Furthermore, if $F=\sum_{\xi\in \hat{U}_d}\alpha_\xi\cdot \xi$, then
    $\alpha_\xi$ is equal to $\alpha_{\tilde{\xi}}$ where $\tilde{\xi}$ is
    the unique lift of $\xi$ to $U_m$ that factors through $U_d$. In particular,
    $F$ is still even and has zero trivial component. Thus, Lemma~\ref{mainlemma}(2)
    tells us that $F$ lies in $V_d$. Since $d$ and $\ff_\chi$ share prime factors,
    Proposition~\ref{eulervanish} allows to conclude that if $\xi\in \hat{U}_d$ such 
    that $\tilde{\xi}=\chi$, then $\xi\in W_d$. Thus, $\alpha_\chi = \alpha_\xi = 0$, which
    concludes our proof.
\end{proof}

\begin{remark}
    In the proof, our choice of $d$ with respect to a Dirichlet character $\chi$ was
    chosen so that $\gcd(m/d,d)=1$ and $d$ shares the same prime factors as $\ff_\chi$.
    The reason we need to choose such a $d$ instead of using the conductor $\ff_\chi$
    is to preserve the even-ness of our projection. If $\gcd(\ff_\chi,m/\ff_\chi)\neq 1$,
    then there exist characters defined modulo $\ff_\chi$ which are not even, but whose
    lift to $U_m$ are even. 
\end{remark}

The theorem above gives rise to the following result which will allow us to fully solve the geometric equidistribution problem in the next section.

\begin{cor}\label{cor:more-general-result}
    Let $\pm c_1,\ldots,\pm c_n,\pm d_1,\ldots,\pm d_n\in U_m$ be two sequences. If 
    \[
    \prod_{i=1}^n (x^{c_i}-1)(x^{-c_i}-1)\equiv \prod_{i=1}^n (x^{d_i}-1) 
    (x^{-d_i}-1)\pmod*{x^m-1},
    \]
    then these two sequences are congruent modulo $m$.
\end{cor}

\begin{proof}
    Let $f,g:U_m\to \Z$ be the functions 
    \begin{align*}
        f(x)&= \#\{1\le i\le n\mid x=\pm c_i\}\\
        g(x)&= \#\{1\le i\le n\mid x=\pm d_i\}.
    \end{align*}
    We see that $f$ and $g$ satisfy the conditions of Theorem~\ref{polyequality},
    so there exists a permutation $\sigma\in S_n$ such that 
    $\pm c_i \equiv \pm d_{\sigma(i)}\pmod*{m}$ for all $1\le i\le n$.
\end{proof}

\section{The geometric equidistribution problem revisited}\label{sec:equidist-revisited}
We are now equipped to characterize when $\gapset(S)$ is equidistributed modulo $m$ for the geometric numerical semigroup $S=\langle \Geom(a,b,;k)\rangle$ for arbitrary $k$. Our strategy will be to first use the structure of the geometric series to find some necessary relations between involving $a,b,k$ and $m$ in Section~\ref{subsec:reduction-to-cases}. Then, we will apply these cases to Eq.~\eqref{eqn:4} to deduce finer relations which are sufficient for equidistribution.

\subsection{Reduction to a few cases}\label{subsec:reduction-to-cases}
As usual, let $a,b,k,m\in\N$ with $a,b>1$, $\gcd(a,b)=1$, $g_i=a^{k-i}b^i$, and $S=\langle \Geom(a,b;k)\rangle=\langle g_0,\dots,g_k\rangle$.

By Proposition~\ref{prop:goal-reframed}, $\gapset(S)$ is equidistributed modulo $m$ if and only if $\gcd(ab,m)=1$ and Eq.~\eqref{eqn:4} holds. Let $r\in\N$ with $r\equiv a^{-1}b\pmod*{m}$. (Such an $r$ necessarily exists because $\gcd(ab,m)=1$.) Then $g_i=g_0r^i$ for all $i$, and we may rewrite Eq.~\eqref{eqn:4} as
\begin{equation}\label{eqn:5}
    (x-1)\prod\limits_{i=1}^k \left(x^{ag_0r^i}-1\right) \equiv \prod\limits_{i=0}^k \left(x^{g_0r^i}-1\right)\pmod*{(x^m-1)}.
\end{equation}

It is tempting to apply Corollary~\ref{cor:more-general-result} to this congruence to conclude that the sequences 
\[
    \pm 1,\, \pm ag_0r,\, \pm ag_0r^2,\, \dots,\, \pm ag_or^k
    \quad\text{ and } \quad 
    \pm g_0r,\, \pm g_0r^2,\, \dots,\, \pm g_0r^k
\]
are congruent modulo $m$, from which we learn that $a$ or $-a$ is a power of $r$ modulo $m$. However, we will take advantage of the geometric structure first. Note that Eq.~\eqref{eqn:5} holds precisely when 
\begin{equation}\label{eqn:5-C_m(x)}
    (x-1)\prod\limits_{i=1}^k \left(x^{ag_0r^i}-1\right) \equiv \prod\limits_{i=0}^k \left(x^{g_0r^i}-1\right)\pmod*{C_m(x)}
\end{equation}
holds. Since $(x^u-1)$ is a unit modulo $C_m(x)$ for any $u\in U_m$, we have 
\begin{equation}\label{eqn:6}
    \prod\limits_{i=1}^k \left(\dfrac{x^{ag_0r^i}-1}{x^{g_0r^i}-1} \right) \equiv \dfrac{x^{g_0}-1}{x-1} \pmod*{C_m(x)}.
\end{equation}
Since $\gcd(r,m)=1$, we may replace $x$ with $x^r$ to obtain
\begin{equation}\label{eqn:7}
    \prod\limits_{i=2}^{k+1} \left(\dfrac{x^{ag_0r^i}-1}{x^{g_0r^i}-1} \right) \equiv \dfrac{x^{g_0r}-1}{x^r-1} \pmod*{C_m(x)}.
\end{equation}
Combining Eqs.~\eqref{eqn:6} and \eqref{eqn:7}, we find 
\begin{equation}
    \dfrac{(x^{ag_0r^{k+1}}-1)}{(x^{g_0r^{k+1}}-1)}\cdot \dfrac{(x^{g_0r}-1)}{(x^{ag_0r}-1)} \equiv \dfrac{(x^{g_0r}-1)}{(x^r-1)}\cdot \dfrac{(x-1)}{(x^{g_0}-1)} \pmod*{C_m(x)}.
\end{equation}

Canceling $(x^{g_0r}-1)$ on each side, we finally obtain 
\begin{equation}\label{eqn:8}
    (x^r-1)(x^{g_0}-1)(x^{ag_0r^{k+1}}-1)\equiv (x-1)(x^{ag_0r}-1)(x^{g_0r^{k+1}}-1)\pmod*{C_m(x)}.
\end{equation}
Hence,
\begin{equation}\label{eqn:9}
    (x^r-1)(x^{g_0}-1)(x^{ag_0r^{k+1}}-1)\equiv (x-1)(x^{ag_0r}-1)(x^{g_0r^{k+1}}-1)\pmod*{x^m-1}.
\end{equation}
The exponent sequences $r,g_0,ag_0r^{k+1}$ and $1,ag_0r,g_0r^{k+1}$ consist of terms which are all units modulo $m$, so by Proposition~\ref{prop:prods-match-exponents-match}, these three-term sequences are congruent modulo $m$. Fixing the orders of these sequences as $r,g_0,ag_0r^{k+1}$ and $1,ag_0r,g_0r^{k+1}$, we can index the six possible congruences between the sequences by the symmetric group $S_3$. After some simple reductions, each congruence simplifies as follows:
\begin{enumerate}
    \item If $\sigma=(1)$, then $r\equiv a\equiv 1\pmod*{m}$.
    \item If $\sigma=(12)$, then $a\equiv 1\pmod*{m}$ and $r$ is unrestricted.
    \item If $\sigma=(23)$, then $r\equiv1\pmod*{m}$ and $a$ is unrestricted.
    \item If $\sigma=(13)$, then $ar\equiv 1\pmod*{m}$ (or, equivalently, $b\equiv 1\pmod*{m}$).
    \item If $\sigma=(123)$, then $a^{k+1}\equiv r^{k+1}\equiv 1\pmod*{m}$ (and $b^{k+1}\equiv1\pmod*{m}$ as well).
    \item If $\sigma=(132)$, then $a^k\equiv r^k\equiv 1\pmod*{m}$ (and $b^{k}\equiv1\pmod*{m}$ as well).
\end{enumerate}

We summarize these results with four cases.
\begin{prop}\label{prop:four-cases}
    For $a,b,k\in\N$ with $\gcd(a,b)=1$, let $S=\langle\Geom(a,b;k)\rangle$. For $m\in\N$, suppose $\gapset(S)$ is equidistributed modulo $m$. Then $\gcd(ab,m)=1$ and (at least) one of the following cases holds:
    \begin{enumerate}
        \item $a\equiv 1\pmod*{m}$ or $b\equiv 1\pmod*{m}$;
        \item $r\equiv 1\pmod*{m}$ (i.e., $a\equiv b\pmod*{m}$);
        \item $r\not\equiv1\pmod*{m}$ and $a^k\equiv b^k\equiv r^k\equiv 1\pmod*{m}$; or
        \item $r\not\equiv1\pmod*{m}$ and $a^{k+1}\equiv b^{k+1}\equiv r^{k+1}\equiv 1\pmod*{m}$.
    \end{enumerate}
\end{prop}
Note that we may assume $r\not\equiv1\pmod*{m}$ in the third and fourth cases because the case where $r\equiv1\pmod*{m}$ is considered in the second case. Additionally, since $ar\equiv b\pmod*{m}$, if any two of $a^t, b^t, r^t$ are 1 modulo $m$ for some $t$, then the third is 1 modulo $m$ as well.

We can plug each of these cases into Eq.~\eqref{eqn:5} to find any additional conditions.

\textbf{Case 1.} If $a\equiv1\pmod*{m}$ or $b\equiv1\pmod*{m}$, then Eq.~\eqref{eqn:5} holds.

\textbf{Case 2.} If $r\equiv1\pmod*{m}$, then by Eq.~\eqref{eqn:5}, $\gapset(S)$ is equidistributed modulo $m$ if and only if 
\begin{equation}\label{eqn:cond2}
    \left(x-1\right)\left(x^{a^{k+1}}-1\right)^k\equiv \left(x^{a^k}-1\right)^{k+1}\pmod*{(x^m-1)}.
\end{equation}

\textbf{Case 3.} If $r\not\equiv1\pmod*{m}$ and $a^k\equiv r^k\equiv 1\pmod*{m}$, then by Eq.~\eqref{eqn:5}, $\gapset(S)$ is equidistributed modulo $m$ if and only if 
\begin{equation}\label{eqn:cond3}
    \prod\limits_{i=0}^{k-1} \left(x^{ar^i}-1\right) 
    \equiv 
    \prod\limits_{i=0}^{k-1} \left(x^{r^i}-1\right)
    \pmod*{(x^m-1)}.
\end{equation}

\textbf{Case 4.} If $r\not\equiv1\pmod*{m}$ and  $a^{k+1}\equiv r^{k+1}\equiv 1\pmod*{m}$, then by replacing $x$ with $x^a$ in Eq.~\eqref{eqn:5}, $\gapset(S)$ is equidistributed modulo $m$ if and only if 
\begin{equation}\label{eqn:cond4}
    \prod\limits_{i=0}^k \left(x^{ar^i}-1\right) 
    \equiv 
    \prod\limits_{i=0}^k \left(x^{r^i}-1\right)
    \pmod*{(x^m-1)}.
\end{equation}

There is nothing more to do with Case 1. For Cases 2, 3, and 4, we need to determine additional conditions on $a$, $r$, and $m$ which imply that $\gapset(S)$ is equidistributed modulo $m$. Our primary tool here is Corollary~\ref{cor:more-general-result}. In what follows, let $R$ be the subgroup of units modulo $m$ generated (multiplicatively) by $r$.

\subsubsection{Case 2}
\begin{lem}
    Suppose $r\equiv1\pmod*{m}$ and $\gcd(a,m)=1$. Then \begin{equation}\label{eqn:cond2-C_m(x)}
        (x-1)\left(x^{a^{k+1}}-1\right)^k
        \equiv 
        \left(x^{a^k}-1\right)^{k+1}
        \pmod*{C_m(x)}.
    \end{equation} 
    if and only if one of the following occurs:
    \begin{itemize}
        \item $a\equiv1\pmod*{m}$; or
        \item $a\equiv -1\pmod*{m}$, $k$ is even, and $m\mid k$; or
        \item $a\equiv -1\pmod*{m}$, $k$ is odd, and $m\mid(k+1)$.
    \end{itemize}
\end{lem}
\begin{proof}
    ($\implies$) Suppose Eq.~\eqref{eqn:cond2-C_m(x)} holds. Since the exponents are all units modulo $m$, we apply Corollary~\ref{cor:more-general-result} to conclude that the sequences 
    \[
        \pm 1,\, \pm a^{k+1},\, \pm a^{k+1},\, \dots,\, \pm a^{k+1}
        \quad\text{ and }\quad 
        \pm a^{k},\, \pm a^{k},\, \pm a^{k},\, \dots,\, \pm a^{k}
    \] 
    are congruent modulo $m$. Thus, $a\equiv \pm 1\pmod*{m}$. If $a\equiv1\pmod*{m}$, then we are done. Otherwise, if $a\equiv -1\pmod*{m}$, then $a^k\equiv (-1)^k\pmod*{m}$. We consider the parity of $k$.

    If $k$ is even, then $a^k\equiv 1\pmod*{m}$ and Eq.~\eqref{eqn:cond2-C_m(x)} becomes $(x-1)(x^{-1}-1)^k\equiv (x^{1}-1)^{k+1}\pmod*{C_m(x)}$. We equivalently have
    \[
        (-x^{-1})^k (x-1)^{k+1}\equiv (x-1)^{k+1}\pmod*{C_m(x)},
    \]
    and thus $x^{-k}\equiv1\pmod*{C_m(x)}$. Equivalently, $x^k\equiv1\pmod*{C_m(x)}$, which implies $C_m(x)\mid (x^k-1)$, and thus $m\mid k$.
    
    If $k$ is odd, then $a^k\equiv -1\pmod*{m}$. Eq.~\eqref{eqn:cond2-C_m(x)} becomes $(x-1)^{k+1}\equiv (x^{-1}-1)^{k+1}\pmod*{C_m(x)}$. Equivalently, 
    \[
        (x-1)^{k+1}\equiv (-x^{-1})^{k+1}(x-1)^{k+1}\pmod*{C_m(x)},
    \]
    and thus $x^{-k-1}\equiv 1\pmod*{C_m(x)}$. Equivalently, $x^{k+1}\equiv1\pmod*{C_m(x)}$, which implies $C_m(x)\mid (x^{k+1}-1)$, and thus $m\mid k+1$.

    ($\impliedby$) Given $a\equiv 1\pmod*{m}$, both sides of Eq.~\eqref{eqn:cond2-C_m(x)} are congruent to $(x-1)^{k+1}$ and hence congruent to each other modulo $C_m(x)$. If $a\equiv -1\pmod*{m}$ and $k$ is even, the right side of Eq.~\eqref{eqn:cond2-C_m(x)} is congruent to $(x-1)^{k+1}$ whereas the left side is congruent to $(x-1)(x^{-1}-1)^k=(x-1)^{k+1}\cdot x^{-k}$. Since $m\mid k$, $x^{-k}\equiv 1\pmod*{C_m(x)}$. The same analysis applies to the case where $k$ is odd and $m$ divides $k+1$.
\end{proof}

\subsubsection{Cases 3 and 4}
\begin{lem}
    Suppose $a^k\equiv r^k\equiv 1\pmod*{m}$ and $r\not\equiv1\pmod*{m}$. Then 
    \begin{equation}\label{eqn:cond3-C_m(x)}
        \prod\limits_{i=0}^{k-1}\left(x^{ar^i}-1\right)
        \equiv
        \prod\limits_{i=0}^{k-1}\left(x^{r^i}-1\right)
        \pmod*{C_m(x)},
    \end{equation}
    if and only if one of the following occurs:
    \begin{itemize}
        \item $a\in R$; or
        \item $-a\in R$, $k$ is even, and $m\mid(1+r+\cdots+r^{k-1})$.
    \end{itemize}
\end{lem}
\begin{proof}
    ($\implies$) By Corollary~\ref{cor:more-general-result}, the sequences 
    \[
        \pm a,\, \pm ar,\, \dots,\, \pm ar^{k-1}
        \quad\text{ and }\quad
        \pm 1,\, \pm r,\, \dots,\, \pm r^{k-1}
    \]
    are congruent modulo $m$. In particular, $\pm a\in R$. If $a\in R$, we are done. If $-a\in R$, then via Eq.~\eqref{eqn:cond3-C_m(x)} we re-index and find
    \[
        \prod\limits_{i=0}^{k-1}\left(x^{-r^i}-1\right)
        \equiv
        \prod\limits_{i=0}^{k-1}\left(x^{r^i}-1\right)
        \pmod*{C_m(x)}.
    \] 
    Equivalently, 
    \[
        1
        \equiv
        \prod\limits_{i=0}^{k-1} \left(-x^{r^i}\right)
        \pmod*{C_m(x)},
    \]
    and therefore
    $(-1)^k x^{1+r+r^2+\cdots+r^{k-1}}\equiv1\pmod*{C_m(x)}$. This implies that $k$ is even and $m\mid (1+r+\cdots+r^{k-1})$.

    ($\impliedby$) If $a\in R$ and $r^k\equiv1\pmod*{m}$, then the sequences $a,ar,ar^2,\dots,ar^{k-1}$ and $1,r,r^2,\dots,r^{k-1}$ are congruent modulo $m$ and hence Eq.~\eqref{eqn:cond3-C_m(x)} holds. 

    If $-a\in R$ with $k$ even and $m\mid\left(1+r+r^2+\cdots+r^{k-1}\right)$, then $\prod_{i=0}^{k-1}(-x^{r^i})\equiv1\pmod*{C_m(x)}$ and thus Eq.~\eqref{eqn:cond3-C_m(x)} holds.
\end{proof}

\begin{remark}
    One may wonder if, in the second case of the above corollary, the condition that $m\mid(1+r+\cdots+r^{k-1})$ is really necessary or if it follows from the other conditions. After all, if $\gapset(S)$ is equidistributed modulo $m$, then $g(S)\equiv0\pmod*{m}$, where 
    \[
        g(S) = \frac{(a-1)(a^{k-1}b+a^{k-2}b^2+\cdots+b^k)-(a^k-1)}{2}
    \]
    by Prop.~\ref{prop:genus}. Since $a^k\equiv1\pmod*{m}$ and $ra\equiv b\pmod*{m}$, we remove $(a^k-1)$ from $g(S)$ and factor a bit to conclude $(a-1)r(1+r+\cdots+r^{k-1})\equiv0\pmod*{m}$. If $(a-1)$ is a unit modulo $m$, then it follows that $1+r+\cdots+r^{k-1}\equiv0\pmod*{m}$. However, this need not always occur!

    For example, let $m=15$, $a=11$, $b=14$, and $k=2$. Then $r\equiv4\not\equiv1\pmod*{15}$. We find $a^k\equiv 11^2\equiv1\pmod*{m}$ and $r^k\equiv 4^2\equiv1\pmod*{m}$. Finally, $-a\equiv-11\equiv4\equiv r\pmod*{m}$< so $-a\in R$. And $k=2$ is even. However, $1+r+\cdots+r^{k-1}=1+r\equiv5\pmod*{5}$, so $m\nmid (1+r+\cdots+r^{k-1})$.
\end{remark}

Finally, for Case 4 we just take our result from Case 3, replacing $k$ with $k+1$. 
\begin{cor}
    Suppose $a^{k+1}\equiv r^{k+1}\equiv 1\pmod*{m}$ and $r\not\equiv1\pmod*{m}$. Let $R$ be the subset of $U_m$ generated by $r$. Then 
    \begin{equation}\label{eqn:cond4-C_m(x)}
        \prod\limits_{i=0}^{k}\left(x^{ar^i}-1\right)
        \equiv
        \prod\limits_{i=0}^{k}\left(x^{r^i}-1\right)
        \pmod*{C_m(x)},
    \end{equation}
    if and only if one of the following occurs:
    \begin{itemize}
        \item $a\in R$; or
        \item $-a\in R$, $k$ is odd, and $m\mid(1+r+\cdots+r^{k})$.
    \end{itemize}
\end{cor}

\subsection{Putting it all together}
We can now write our general result.

\begin{thm}\label{thm:main-result}
    For $a,b,k,m\in\N$ with $a,b>1$ and $\gcd(a,b)=1$, and for the geometric numerical semigroup $S=\langle \Geom(a,b;k) \rangle$, the set of gaps of $S$ is equidistributed modulo $m$ if and only if $\gcd(ab,m)=1$ and, for $r$ an integer with $r\equiv a^{-1}b\pmod*{m}$ and $R$ the subgroup of $U_m$ generated (multiplicatively) by $r$, at least one of the following conditions holds:
    \begin{enumerate}
        \item $a\equiv 1\pmod*{m}$ or $b\equiv 1\pmod*{m}$;
        \item $a\equiv b\equiv -1\pmod*{m}$, $k$ is even, and $m\mid k$;
        \item $a\equiv b\equiv -1\pmod*{m}$, $k$ is odd, and $m\mid(k+1)$;
        \item $a\not\equiv b\pmod*{m}$, $a^k\equiv b^k\equiv1\pmod*{m}$, and $a\in R$;
        \item $a\not\equiv b\pmod*{m}$, $a^k\equiv b^k\equiv1\pmod*{m}$, $-a\in R$, $k$ is even, and $m\mid (1+r+\cdots+r^{k-1})$;
        %\item $a^{k+1}\equiv r^{k+1}\equiv 1\pmod*{m}$ and $a\in R$;
        \item $a\not\equiv b\pmod*{m}$, $a^{k+1}\equiv b^{k+1}\equiv1\pmod*{m}$, and $a\in R$;
        %\item $a^{k+1}\equiv r^{k+1}\equiv 1\pmod*{m}$, $-a\in R$, $k$ is odd, and $m\mid (1+r+\cdots+r^{k})$.
        \item $a\not\equiv b\pmod*{m}$, $a^{k+1}\equiv b^{k+1}\pmod*{m}$, $-a\in R$, $k$ is odd, and $m\mid (1+r+\cdots+r^{k})$.
    \end{enumerate}
    %In the conditions above, $r$ is an integer such that $r\equiv a^{-1}b\pmod*{m}$, and $R$ is the subgroup of $U_m$ generated by $r$.
\end{thm}

We encourage the reader to confirm that when $k=1$ or $k=2$, we recover the previously known results for $S=\langle a,b\rangle$ (Proposition~\ref{prop:shor-equidistribution}) and $S=\langle a^2,ab,b^2\rangle$ (Proposition~\ref{prop:result-k=2}).

\begin{remark}
    While conditions 1 through 3 are visually symmetric with respect to $a$ and $b$, it happens that all seven conditions are symmetric with respect to $a$ and $b$. To see this, let $r\equiv a^{-1}b\pmod*{m}$ and $r'\equiv b^{-1}a\pmod*{m}$, which we observe are inverses of one another so that the set of powers of $r$ and the powers of $r'$ in $U_m$ are equal. Thus, it is sufficient to show the following three facts:
    \begin{itemize}
        \item $a\in R$ if and only if $b\in R$,
        \item $-a\in R$ if and only if $-b\in R$,
        \item for all $k\in\N$, $m\mid (1+r+\cdots+r^{k})$ if and only if $m\mid (1+r'+\cdots+(r')^{k})$.
    \end{itemize}
    The first fact is straightforward: $R$ is closed under multiplication by $r$, so if $a\in R$, then $ra\equiv a^{-1}ba\equiv b$, so $b\in R$. The same argument holds for the second fact. For the last fact, we see that if $m\mid (1+r+\cdots+r^k)$, then $m$ must divide $(r')^k(1+r+\cdots+r^k)=(r')^k+\cdots+r'+1$ since $r'\equiv r^{-1}\pmod*{m}$. By the symmetry of this argument, the converse is also true.
\end{remark}

\begin{example}
    To emphasize the necessity of all 7 conditions in the statement of Theorem~\ref{thm:main-result}, we provide for each condition a family of examples. Each family exclusively satisfies the respective condition and we include arguments for why no other of the conditions is satisfied.
    \begin{enumerate}
        \item $a=13$, $b=7$, $m=3$, and $k\in \N$. Since $a\equiv b\equiv 1\pmod*{3}$, no other conditions are satisfied.
        \item $a=5$, $b=11$, $m=3$, and $6\mid k$. Since $a\equiv b\equiv -1\pmod*{3}$ and $k$ is even, no other conditions are satisfied.
        \item $a=5$, $b=11$, $m=3$, and $6\mid (k+1)$. Since $a\equiv b\equiv -1\pmod*{3}$ and $k$ is odd, no other conditions are satisfied.
        \item $a=9$, $b=4$, $m=7$, and $3\mid k$. We have $r\equiv 2\pmod*{7}$ which generates the squares modulo $7$. Since $-2$ is not a square modulo $7$, conditions 5 and 7 are not satisfied. Both $a$ and $b$ have order $3$ modulo $7$, so if $3\mid k$, then $a^{k+1}\equiv a\pmod*{7}$ and $b^{k+1}\equiv b\pmod*{7}$, neither of which are $1\pmod*{7}$. This means condition 6 is not satisfied. Since $a$ and $b$ are each not congruent to $\pm 1\pmod*{7}$, conditions 1 through 3 are also not satisfied.
        \item $a=3$, $b=13$, $m=7$, and $6\mid k$. We have $r\equiv 2\pmod*{7}$, so as in the previous example, $R$ is the set of squares modulo $7$. With $a$ being a generator of $U_7$, $a\notin R$, so conditions 4 and 6 are not satisfied. Condition 7 is clearly not satisfied since $k$ is even in addition to the fact that $1+r+\cdots+r^{k}\equiv 2\pmod*{7}$. Both $a$ and $b$ are not congruent to $\pm 1\pmod*{7}$, so conditions 1 through 3 are also not satisfied.
        \item $a=9$, $b=13$, $m=7$, and $3\mid (k+1)$. The exact argument for example 4 with appropriate adjustments show this example only satisfies condition 6.
        \item $a=3$, $b=13$, $m=7$, and $6\mid (k+1)$. The exact argument for example 5 with appropriate adjustments show this example only satisfies condition 7.
    \end{enumerate}
\end{example}

\section{Proof that \texorpdfstring{$ab$}{\textit{ab}} and \texorpdfstring{$m$}{\textit{m}} must be coprime}\label{sec:gcd(ab,m)=1}
In this section, we provide a proof of Proposition~\ref{prop:gcd(ab,m)=1}, which states that if the set of gaps of a geometric numerical semigroup is equidistributed modulo $m$, then $m$ is coprime to the generating elements. Our approach is to let $\gcd(a,m)=d$ and assume the set of gaps of $S$ is equidistributed modulo $m$. This implies the set of gaps is equidistributed modulo $d$. Thus, the number of gaps in a particular congruence class modulo $d$ is equal to the genus of $S$ divided by $d$. We will explicitly compute the number of gaps that are 0 modulo $d$, set that equal to $g(S)/d$, and conclude that $d=1$. Since $a$ and $b$ are interchangeable, we can similarly show that $\gcd(b,m)=1$ and arrive at our desired result.

As usual, suppose $\gcd(a,b)=1$, let $g_i=a^{k-i}b^i$ for $0\le i\le k$, and consider $S=\langle\Geom(a,b;k)\rangle=\langle g_0,\dots,g_k\rangle$, a geometric numerical semigroup. We will use $f(x)=nx^n$ and $g_0\in S$ with Proposition~\ref{prop:tuenter-generalization} to obtain
\[
    \sum\limits_{n\in \gapset(S)}\left((n+g_0)x^{n+g_0}-nx^n\right) 
    = \sum\limits_{n\in\Ap(S;g_0)}nx^n 
    - \sum\limits_{n=0}^{g_0-1}nx^n.
\]
Given our explicit description of $\Ap(S;g_0)$ from Proposition~\ref{prop:apery-set}, we have that
\begin{equation}\label{eqn:non-coprime-tuenter}
    \sum\limits_{n\in \gapset(S)}\left((n+g_0)x^{n+g_0}-nx^n\right) 
    = \sum\limits_{n_1=0}^{a-1}\cdots\sum\limits_{n_k=0}^{a-1}(n_1g_1+\cdots+n_kg_k)x^{n_1g_1+\cdots+n_kg_k} 
    - \sum\limits_{n=0}^{g_0-1}nx^n.
\end{equation}

Let $d=\gcd(a,m)$. Then $d$ divides $a$ as well as $g_i$ for $0\le i\le k-1$. We will work modulo $(x^d-1)$. Our goal is to write the polynomial $\sum_{n\in \gapset(S)}x^n$ reduced modulo $(x^d-1)$ as 
\[\sum\limits_{n\in \gapset(S)}x^n = \sum\limits_{i=0}^{d-1}c_ix^i\pmod*{(x^d-1)}\]
for some coefficients $c_0,\dots,c_{d-1}$ and then to explicitly calculate $c_0$, the number of gaps of $S$ that are congruent to 0 modulo $d$.

Consider the polynomials $f_1(x)$, $f_2(x)$, $f(x)$ defined as follows:
\begin{itemize}
    \item $f_1(x)=\sum\limits_{n_1=0}^{a-1}\sum\limits_{n_2=0}^{a-1}\cdots\sum\limits_{n_k=0}^{a-1}(n_1g_1+n_2g_2+\cdots+n_kg_k)x^{n_1g_1+n_2g_2+\cdots+n_kg_k}$,
    \item $f_2(x)=\sum\limits_{n=0}^{g_0-1}nx^n$,
    \item $f(x)=f_1(x)-f_2(x)$.
\end{itemize}
Note that $f(x)$ is the right side of Eq.~\eqref{eqn:non-coprime-tuenter}. We wish to reduce $f(x)$ modulo $(x^d-1)$ and compute its constant term, which will then allow us to calculate $c_0$. We will do so by computing the constant terms of the reductions of $f_1(x)$ and $f_2(x)$ modulo $(x^d-1)$ in Lemmas~\ref{lem:f1} and \ref{lem:f2} below.

\begin{lem}\label{lem:f1}
    Suppose $d$ divides $a$. Then $f_1(x)\equiv \sum\limits_{i=0}^{d-1}c_{1,i} x^i\pmod*{(x^d-1)}$ for some integers $c_{1,0},\dots,c_{1,d-1}$. We have 
    \[c_{1,0} = a^k\cdot \frac{b^k(ab-a+ad-bd)-a^k(a-1)b}{2(b-a)d}.\]
\end{lem}
\begin{proof}
    To begin, since $d$ divides $a$ we have $x^a\equiv1\pmod*{(x^d-1)}$. Hence, 
    \[x^{n_1g_1+\cdots+n_kg_k}\equiv x^{n_kg_k}\pmod*{(x^d-1)}.\]
    
    Next, we focus on the first $k-1$ summations in $f_1(x)$. We have that 
    \[
        \sum\limits_{n_1=0}^{a-1}\cdots\sum\limits_{n_{k-1}=0}^{a-1}n_1g_1
        = \sum\limits_{n_1=0}^{a-1}n_1g_1\sum\limits_{n_2=0}^{a-1}1\sum\limits_{n_3=0}^{a-1}1\cdots\sum\limits_{n_{k-1}=0}^{a-1}1
        =\frac{a^{k-1}(a-1)g_1}{2}.
    \]
    Similarly, for all $0<i<k$, 
    \[
        \sum\limits_{n_1=0}^{a-1}\cdots\sum\limits_{n_{k-1}=0}^{a-1}n_ig_i
        =\frac{a^{k-1}(a-1)g_i}{2}.
    \]

    We now look at the final summation in $f_1(x)$, which is over the index $n_k$. Since $\gcd(a,b)=1$ and $d$ divides $a$, we have $\gcd(b,d)=1$ and hence $\gcd(b^k,d)=1$. Thus, $b^k$ is a unit modulo $d$, and so $n_kg_k=n_kb^k\equiv0\pmod*{d}$ precisely when $d$ divides $n_k$. Equivalently, $x^{n_kg_k}\equiv1\pmod*{(x^d-1)}$ precisely when $d$ divides $n_k$. Therefore, in this final summation, since we are calculating the constant term of $f_1(x)$ modulo $(x^d-1)$, we need only consider values of $n_k$ among the $a/d$ values $0, d, 2d, \dots, (a/d-1)d$. It follows that the constant term (modulo $(x^d-1)$) of $\sum_{n_k=0}^{a-1}x^{n_kg_k}$ is $a/d$.
    
    For $0<i<k$, the constant term of 
    \[
        \sum\limits_{n_1=0}^{a-1}\cdots\sum\limits_{n_{k-1}=0}^{a-1}\sum\limits_{n_{k}=0}^{a-1}n_ig_ix^{n_kg_k} \text{ is } \frac{a^{k-1}(a-1)g_i}{2}\cdot\frac{a}{d}=\frac{a^k(a-1)g_i}{2d}.
    \]
    For $i=k$, the constant term of 
    \[
        \sum\limits_{n_1=0}^{a-1}\cdots\sum\limits_{n_{k-1}=0}^{a-1}\sum\limits_{n_{k}=0}^{a-1}n_kg_kx^{n_kg_k}
    \]
    is
    \[ 
        a^{k-1}(0g_k+dg_k+2dg_k+\cdots+(a/d-1)dg_k)
        =\frac{a^{k}(a-d)g_k}{2d}.
    \]
    
    We can then sum from $i=1$ to $i=k$ to obtain the constant term of $f_1(x)$:
    \[
        \sum\limits_{i=1}^{k-1}\frac{a^k(a-1)g_i}{2d} + \frac{a^k(a-d)g_k}{2d}
        =\frac{a^k(a-1)}{2d}\frac{ab(b^{k-1}-a^{k-1})}{b-a} + \frac{a^kb^k(a-d)}{2d}.
    \]
    In this calculation, we used the fact that $g_1+g_2+\cdots+g_{k-1}$ is geometric, equal to $ab(b^{k-1}-a^{k-1})/(b-a)$. Algebraic manipulation gets us the stated result.
\end{proof}
\begin{lem}\label{lem:f2}
    Suppose $d$ divides $a$. 
    Then $f_2(x)\equiv \sum\limits_{i=0}^{d-1}c_{2,i} x^i\pmod*{(x^d-1)}$ for some integers $c_{2,0},\dots,c_{2,d-1}$. We have 
    \[c_{2,0} = a^k\cdot \frac{a^k-d}{2d}.\]
\end{lem}
\begin{proof}
    Since we are working modulo $(x^d-1)$, the monomials of $f_2(x)$ that reduce to $x^0$ are exactly those with exponents divisible by $d$. Since the summation in $f_2(x)$ goes from $n=0$ to $g_0-1=a^k-1$, we get the following finite arithmetic series with common difference $d$:
    \[c_{2,0}=0+d+2d+\cdots+(a^k/d - 1)d=\frac{(a^k-d)a^k}{2d}.\qedhere\]
\end{proof}
\begin{lem}\label{lem:num_gaps_0_mod_d}
    Suppose $d$ divides $a$. The number of gaps of $S=\langle \Geom(a,b;k)\rangle$ that are divisible by $d$ is
    \[\frac{(a-1)b^{k+1}-(b-1)a^{k+1}+(b-a)}{2(b-a)d}+\frac{b^k(1-d)}{2d}.\]
\end{lem}
\begin{proof}
    We begin by considering Eq.~\eqref{eqn:non-coprime-tuenter} modulo $(x^d-1)$. As we noted earlier, we have $x^a\equiv 1\pmod*{(x^d-1)}$ and $x^{g_i}\equiv1\pmod*{(x^d-1)}$ for $i=1,2,\dots,k-1$. We obtain
    \[
        a^k\sum\limits_{n\in \gapset(S)}x^n 
        \equiv \sum\limits_{n_1=0}^{a-1}\cdots\sum\limits_{n_k=0}^{a-1}(n_1g_1+\cdots+n_kg_k)x^{n_kg_k} 
        - \sum\limits_{n=0}^{a^k-1}nx^n \pmod*{(x^d-1)}.
    \]
    Note that $\sum_{n\in \gapset(S)}x^n$ is congruent to a unique polynomial of the form $c_0+c_1x+\cdots+c_{d-1}x^{d-1}$ modulo $(x^d-1)$, and each coefficient $c_r$ is equal to the number of gaps of $S$ that are congruent to $r$ modulo $d$. For this proof, we are interested in computing $c_0$.
    
    For $f_1(x)$ and $f_2(x)$ as defined in Lemmas~\ref{lem:f1} and \ref{lem:f2}, we have $a^k\sum_{n\in \gapset(S)}x^n = f_1(x)-f_2(x)$. Hence, $a^kc_0=c_{1,0}-c_{2,0}$, from which we find
    \[c_0 = \frac{b^k(ab-a+ad-bd)-a^k(a-1)b}{2(b-a)d} - \frac{(a^k-d)}{2d}.\]
    Via algebraic manipulation, we find that 
    \[c_0 = \frac{(a-1)b^{k+1}-(b-1)a^{k+1}+(b-a)d}{2(b-a)d}+\frac{b^k(1-d)}{2d}.\]
    This is the number of gaps of $S$ that are congruent to 0 modulo $d$.
\end{proof}

We are now able to prove that, for $S=\langle\Geom(a,b;k)\rangle$, if $\gapset(S)$ is equidistributed modulo $m$, then $\gcd(ab,m)=1$.

\begin{proof}[Proof of Proposition~\ref{prop:gcd(ab,m)=1}]
    Suppose $\gapset(S)$ is equidistributed modulo $m$. %If $k=1$, then $S=\langle a,b\rangle$. By Proposition~\ref{prop:shor-equidistribution}, we have $\gcd(ab,m)=1$.
    %Now, suppose $k\ge2$. 
    Let $\gcd(a,m)=d$. Since $d$ divides $m$, the set of gaps of $S$ is also equidistributed modulo $d$. This means there are $g(S)/d$ gaps in each congruence class modulo $d$. With the formula for $g(S)$ in Prop.~\ref{prop:genus}, we obtain
    \[g(S)/d = \frac{(a-1)b^{k+1}-(b-1)a^{k+1}+(b-a)}{2(b-a)d}.\]

    By Lemma~\ref{lem:num_gaps_0_mod_d}, since $d$ divides $a$, the number of gaps that are congruent to 0 modulo $d$ is \[\frac{(a-1)b^{k+1}-(b-1)a^{k+1}+(b-a)d}{2(b-a)d} + \frac{b^k(1-d)}{2d}.\] 

    These quantities are equal. Canceling terms in common, we have 
    \[\frac{(b-a)}{2(b-a)d} = \frac{(b-a)d}{2(b-a)d}+\frac{b^k(1-d)}{2d}.\]
    Equivalently, $(b^k-1)(1-d)=0$. Since $b^k\ge2$, we must have $d=1$. Thus, $\gcd(a,m)=1$ as desired.

    Finally, note that we have made no assumptions about $a$ and $b$ relative to each other in this section. We may therefore swap their roles and follow the same argument to conclude that $\gcd(b,m)=1$.
\end{proof}

\bibliographystyle{plainnat}
%\bibliography{refs}

\end{document}